\documentclass[twoside,12pt,a4paper,mathscr]{article}
\usepackage{amsmath,amssymb,amsthm,euscript,amscd,a4wide}
\usepackage[all]{xypic}
\CompileMatrices
\newtheorem{thm}{Theorem}[section]

\newtheorem{corol}[thm]{Corollary}
\newtheorem{lemma}[thm]{Lemma}
\newtheorem{prop}[thm]{Proposition}
\newtheorem{defin}[thm]{Definition}
\newtheorem{defin1}[thm]{Definition (I)}
\newtheorem{defin2}[thm]{Definition (II)}
\newtheorem{rem}[thm]{Remark}
\newtheorem{ex}[thm]{Example}

\newenvironment{remark}{\begin{rem}\rm}
       {\hfill$\vartriangle$\mbox{\hskip3pt}\end{rem}}
\let\cal=\mathcal
\newcommand{\Pc}{{\cal P}}

\newcommand{\C}{{\mathbb C}}

\newcommand{\dbar}{\overline\partial}
\newcommand{\lra}{\longrightarrow}

\newcommand{\Rb}{{\mathbf R}}
\newcommand{\Hom}{\operatorname{Hom}}

\newcommand{\Ext}{\operatorname{Ext}}
\newcommand{\Pic}{\operatorname{Pic}}

\newcommand{\Ker}{\operatorname{Ker}}
\newcommand{\Coker}{\operatorname{Coker}}

\newcommand\ch{\operatorname{ch}}
\newcommand\WIT{\operatorname{WIT}}
\newcommand\IT{\operatorname{IT}}
\newcommand\rk{\operatorname{rk}}

\newcommand\Id{\operatorname{Id}}

\newcommand\tr{\operatorname{tr}}
\newcommand\iso{\xrightarrow{\sim}}
\newcommand{\bysame}{$\raise.2em\hbox to 3em{\hrulefill}$\thinspace, }
\begin{document}
\thispagestyle{empty}
\begin{center}
{\bfseries\Large  Fourier-Mukai and Nahm transforms for holomorphic
triples on elliptic curves}
\par\addvspace{20pt}
{\sc Oscar Garc\'ia-Prada}$^{\P}$,\ {\sc Daniel Hern\'andez Ruip\'erez}$^{\dag}$,\
{\sc Fabio Pioli}$^{\P}$\\{\sc Carlos Tejero Prieto}$^{\S}$
\par\medskip
\P\ Instituto de Matem\'aticas y F\'isica Fundamental,\\
Consejo Superior de Investigaciones Cient\'{\i}ficas,\\
Serrano 113 bis, 28006 Madrid, Spain
\par\medskip
\dag\ Departamento de Matem\'aticas and Instituto Universitario de
F\'{\i}sica Fundamental y Matem\'aticas, Universidad de Salamanca, \\
Plaza de la Merced 1-4, 37008 Salamanca, Spain\\
\par\medskip
\S\  Departamento de Matem\'aticas, Universidad de Salamanca, \\
Plaza de la Merced 1-4, 37008 Salamanca, Spain
\par\medskip
\end{center}
\vfill
\begin{quote} \footnotesize {\sc Abstract.}
We define a Fourier-Mukai transform for a triple
consisting of two holomorphic vector bundles over an elliptic
curve and a homomorphism between them. We
prove that in some cases the transform preserves the natural stability
condition for a triple. We also define a Nahm transform for
solutions to natural gauge-theoretic equations on a triple --- vortices ---
and explore some of its basic properties. Our approach combines direct methods
with dimensional reduction techniques, relating triples over a curve with
vector bundles over the product of the curve with the complex projective line.
\end{quote}
\vfill
\leftline{\hbox to8cm{\hrulefill}}\par
{\footnotesize
\noindent  The authors are members of VBAC (Vector bundles on
algebraic curves), which is partially supported by EAGER (EC FP5
Contract no. HPRN-CT-2000-00099) and by EDGE (EC FP5 Contract no.
HPRN-CT-2000-00101). This research has been partially supported
by the Italian/Spain bilateral programme Azione Integrata, IT203
"Sheaves on Calabi-Yau manifolds and applications to integrable
systems and string theory" and by the research projects BFM2003-00097 of the
spanish DGI and SA118/03 of the ``Junta de Castilla y Le\'on''.
\\
\noindent\emph{E-Mail addresses:} {\tt oscar.garcia-prada@uam.es, ruiperez@usal.es,
fabio.p@imaff.cfmac.csic.es, carlost@usal.es} 
}

\eject\thispagestyle{empty}\mbox{\ \ \ }
\setcounter{page}{1}

\section{Introduction}

The Fourier-Mukai transform, as originally introduced by Mukai for abelian
varieties  \cite{Mu} establishes a duality between the derived categories of
coherent sheaves over an abelian variety and its dual variety. The theory has
been extended to more general varieties, including K3 surfaces, Calabi-Yau
threefolds or elliptic fibrations. In particular, it  is a very powerful tool
in the study of moduli spaces of vector bundles over abelian surfaces and K3
surfaces (see \cite{Mu1, Br1, Mac, HP, Yo} for instance). In the
gauge-theoretic side, the Nahm transform provides a differential geometric
analogue of the Fourier-Mukai transform relating instantons (or monopoles) on
dual manifolds \cite{Na, Hi}. In many cases, whenever it makes sense, both
transforms are compatible in a suitable way.

In this paper we study Fourier-Mukai and Nahm transforms for holomorphic
triples over an elliptic curve  and their corresponding vortex equations.
A triple here consists of two holomorphic vector bundles over
the elliptic curve and a homomorphism between them.
The motivation to study this problem is two-fold. On the one hand the
Nahm transform has been successfully applied to find  instanton and
monopole solutions, which are defined  in real dimensions
4 and 3 respectively. It is then very natural to try to find an analogue
for 2-dimensional vortices. On the other hand, vortices in two dimensions
are equivalent to $SU(2)$-invariant instantons over the product of the
elliptic curve and the Riemann sphere, where the $SU(2)$ action is given
simply by the usual one on the sphere. This suggests a relative 4-dimensional
approach to the problem.  In a related context the Nahm transform
has been successfully applied to study doubly periodic instantons and their
relationship with Hitchin's equations \cite{J1,J2}.

Here is a description of the paper. In  Section \ref{ftransform} we briefly
review the Fourier-Mukai and Nahm transforms for vector bundles over elliptic
curves. We recall the preservation of stability and prove that the constant
central curvature condition for a connection (which on a curve
coincides with the Einstein-Hermitian condition) is preserved.
Although the latter seems
to be of general knowledge, we have not found it in the literature and hence
include it here since it is relevant for our further study for triples.
We follow the approach given in \cite{G}.

In Section \ref{transtrip}, we review first the basic stability theory for
triples. An important feature is that the stability criterium depends on a
real parameter which is typically bounded \cite{BG}. We then  introduce the
Fourier-Mukai transform for triples on elliptic curves and give  two natural
approaches for transforming a triple. The first one is based on the absolute
Fourier-Mukai transform acting on the components of the triple. The second
approach is based on a relative Fourier-Mukai transform combined with a
dimensional reduction procedure. We prove that the Fourier-Mukai transform
preserves stability of triples for ``small'' and ``large'' values of the
stability parameter, providing an isomorphism of moduli spaces. What happens
for other values of the parameter remains to be investigated.  We conclude
this section by applying these results to obtain isomorphisms between moduli
spaces of stable $SU(2)$-equivariant vector bundles.

Finally, in Section \ref{transnahm}, in parallel with Section \ref{transtrip},
we develop the formalism for a relative Nahm transform in the same context. We
apply this formalism to transform a solution to the vortex equations over a
triple, regarded as an $SU(2)$-invariant Einstein-Hermitian connection on a
certain $SU(2)$-equivariant bundle over the product of the curve with the
complex projective line. In general it seems very hard to identify   the
equation satisfied by the Nahm transform of a vortex solution, which one would
expect to be again the vortex equation on the transformed triple. We content
ourselves with analysing in full detail the case of covariantly constant
triples, leaving for a future paper the analysis of the general case. As a
byproduct we prove that polystability of triples may not be preserved by the
Fourier-Mukai transform.

In this paper we work over the field of complex numbers $\mathbb C$.

\section{Fourier-Mukai and Nahm transforms on elliptic cur\-ves}
\label{ftransform}

\subsection{Fourier-Mukai transform}\label{section:FM}

Let $C$ be an elliptic curve and let $\widehat C=\Pic^0 (C)$ be its dual
variety. Although $C$ and $\widehat C$ are isomorphic it will be convenient to
keep a notational distinction between them for the sake of clarity. Over
$C\times \widehat C$ we consider the Poincar\'e bundle $\Pc $ and we denote by
$\pi_C$ and $\pi_{\widehat C}$ the canonical projections onto the factors $C$
and $\widehat C$. As it is customary, $\Pc$ is normalized so that it is trivial
over $\{0\}\times \widehat C$. In  \cite{Mu} Mukai introduced a functor between
the bounded derived categories of coherent sheaves of $C$ and $\widehat C$: $$
\mathcal S \colon
 D(C) \to D(\widehat C).$$ This functor acts as follows
$$\mathcal S (E) = \Rb \pi_{\widehat C,*} (\pi_C^* E \otimes \Pc), $$ where
$E$ is an object of the derived category and $\Rb\pi_{\widehat C,*}$ denotes
the derived functor of $\pi_{\widehat C,*}$.

We need some standard terminology and notation. As usual, we denote by $\mathcal
S^i(E)$ the sheaf defined by the $i$-th cohomology of the
complex $\mathcal S(E)$, that is,
 $$\mathcal S^i(E) = \mathcal{H}^i (\mathcal S
(E))\,. $$ When $E$ is a single sheaf, $S^i(E)$ is the ordinary derived
functor $R^i\pi_{\widehat C,*}(\pi_C^* E\otimes \Pc)$. A sheaf $E$ is said to
be $\WIT_i$ if $\mathcal S^j(E) = 0$ for every $j\neq i$, and $E$ is called
$\IT_i$ if it is $\WIT_i$ and its unique  transform $\mathcal S^i (E)$ is
locally-free. Equivalently $E$ is $\IT_i$ if the cohomology group $H^j (C_\xi,
E \otimes \Pc_\xi) = 0$ vanishes for every $j\neq i$ and every $\xi \in
\widehat C$, where $C_\xi = C\times {\{\xi\}}$ and $\Pc_\xi $ is the
restriction of $\Pc $ to $C_\xi$. In this case, the fibre over $\xi\in\widehat
C$ of the unique Fourier-Mukai transform $S^i(E)$ is canonically isomorphic to
$H^i(C_\xi,E\otimes \Pc_\xi)$. The Fourier-Mukai transform $S^i(E)$ of a
$\WIT_i$ sheaf $E$ will be denoted as usual by $\widehat E$. When there is no
need to specify the index $i$ we shall simply say that a sheaf is $\WIT$ or
$\IT$.

One of the most important features of the functor $S$ is that it admits an
inverse $\hat{\mathcal S}\colon D(\widehat C) \to D( C)$. That is, there are natural
isomorphisms:
\begin{equation*}
\begin{split}
 \hat{\mathcal S} &\circ \mathcal S \simeq \Id_{D(C)}\\
\mathcal S &\circ \hat{\mathcal S}  \simeq \Id_{D(\widehat C)}.
\end{split}
\end{equation*} Moreover $\hat{\mathcal S}$ is explicitly given by
$$\hat{\mathcal S} (F) =  \Rb \pi_{C,*} (\pi_{\widehat C}^* (F) \otimes
\Pc^\vee[1]),$$ where $\Pc^\vee$ is the dual of $\Pc$ and $[1]$ is the shift
operator.

Let us recall the following well-known fact whose proof relies on the
invertibility property of the Fourier-Mukai transform (see \cite{G},
also \cite{Dekk} and \cite{Br1}).

\begin{prop}\label{presstab} If $E$ is a  semistable (stable)
vector bundle of non-zero degree over an elliptic curve $C$, then $E$ is $\IT$
and the transform $\widehat E$ is also semistable (stable). Moreover, $E$ is
$\IT_0$ ($\IT_1$) if and only if $\deg(E)>0$ ($\deg(E)<0$). Finally, if $E$ is
$\IT_i$ with Chern character $\ch(E)=(r,d)$ then $\ch(\widehat E)=((-1)^i
d,(-1)^{i+1}r)=(-1)^i(d,-r)$.
\qed\end{prop}

\begin{remark} If we take into account that any vector bundle $E$ on
an elliptic curve decomposes uniquely into a direct sum of semistable bundles
we conclude that $E$ is $\IT_0$ ($\IT_1$) if and only if all of its components
have positive (negative) degree.
\end{remark}

Recall that on an elliptic curve $C$ the moduli space $\mathcal M_C (r,d)$ of
$S$-equivalence classes of semistable bundles of rank $r$ and degree $d$ is
isomorphic to the symmetric product $S^h C$, where $h=(r,d)$ is the greatest
common divisor of $r$ and $d$. When $(r,d)>1$ there  are no stable bundles in
$\mathcal M_C (r,d)$. When $r$ and $d$ are coprime, all the semistable bundles
are stable and $\mathcal M_C (r,d)$ is isomorphic to $C$ (see \cite{At} and
\cite{Tu} for details, as well as \cite{Br1} and \cite{HP}). The Fourier-Mukai
transform is well-behaved with respect to families of stable bundles and
therefore induces morphisms between moduli spaces. In the case of $\IT_i$
semistable bundles on an elliptic curve, the Fourier-Mukai transform also
preserves $S$-equivalence. More precisely if $E$ is an $\IT_i$ semistable
bundle on $C$, then it is immediate to see that every stable bundle $E_k$ in
the graded object $\mathrm{Gr} (E) = \oplus_k E_k$ with respect to a Jordan-
H\"older filtration is $\IT_i$. From this follows that if $E$ and $E^\prime$
are $S$-equivalent $\IT_i$ bundles, then the transforms $\widehat E$ and
$\widehat E^\prime$ remain $S$-equivalent. Therefore we have.

\begin{corol} Let $\mathcal M_C (r,d)$ be the moduli space of
semistable bundles of rank $r$ and $d\neq 0$. Then, in the $\IT_i$
case, the Fourier-Mukai transform induces an isomorphism between the
moduli spaces
$$\mathcal S \colon \mathcal M_C (r,d)\iso \mathcal M_{\widehat C}
((-1)^id, (-1)^{i+1} r).$$
Therefore the Fourier-Mukai transform gives rise to an isomorphism
between symmetric products of elliptic curves.
\end{corol}

\subsection{Nahm transform}\label{nahm}

We come now to the definition of the Nahm transform in the context of elliptic
curves.

Let $C$ be a complex elliptic curve endowed with a flat metric of unit volume.
The canonical spinor bundle $S=\Lambda^{0,\bullet}\,T^*C$ of $C$ as a spin$^c$
manifold, has a natural splitting $S=S^+\oplus S^-$ where
$$S^+=\Lambda^{0,0}T^*C,\qquad S^-=\Lambda^{0,1}T^*C\,.$$ We denote the
spinorial connection of $S$ by $\nabla_S$.

The dual elliptic curve $\widehat C$ parametrizes the gauge equivalence classes
of Hermitian flat line bundles over $C$. The Poincar\'e bundle $\mathcal{P}$
introduced in Section \ref{section:FM} is endowed with a unitary connection
$\nabla_\mathcal{P}$, such that the restriction  of
$(\mathcal{P},\nabla_\mathcal{P})$ to the slice $C_{\xi}$ is in the
equivalence class defined by $\xi\in\widehat C$. Therefore for every
$\xi\in\widehat C$ we have the Hermitian line bundle
$\mathcal{P}_{\xi}\equiv\mathcal{P}_{\vert{C_{\xi}}}\to C$ endowed with the
flat unitary connection
$\overline\nabla_{\xi}=\nabla_{\vert\mathcal{P}_{\xi}}$.

Let us consider a Hermitian vector bundle $E\to C$ with a unitary connection
$\nabla$. On the vector bundle $E\otimes\mathcal{P}_{\xi}$ we have the
connection $\nabla_{\xi}=\nabla\otimes 1+1\otimes \overline\nabla_{\xi}$.
Therefore we have the family of coupled Dirac operators
$$D_{\xi}\colon\Omega^0(C,S^+\otimes E\otimes \mathcal{P}_{\xi})\to
\Omega^0(C,S^-\otimes E\otimes \mathcal{P}_{\xi}).$$ It follows from the
Atiyah-Singer Theorem  for families that the difference bundle of
the family of Dirac operators $D$ parametrized by $\widehat C$ is a well defined
object $\mathrm{Ind}(D)$ in $K$-theory which is called the index of $D$.
Moreover, if either one of $\{\Ker  D_{\xi}\}$ or $\{\Coker
D_{\xi}\}$ has constant rank, then $\Ker D$ and $\Coker  D$ are
vector bundles over $\widehat C$ and one has that

$$\mathrm{Ind}(D)= [\Ker  D]-[\Coker  D]\in K(\widehat C)\,.$$

\begin{defin}
Let $(E,\nabla)$ be a pair formed by a Hermitian vector bundle $E$ over $C$
and a unitary connection $\nabla$ on $E$. We say that $(E,\nabla)$ is an
$\IT$ (index Theorem) pair if either $\Coker D=0$ or $\Ker D=0$.
In the first case we say that $(E,\nabla)$ is an $\IT_0$-pair whereas in the
second we call it an $\IT_1$-pair. The transformed bundle of an $\IT_i$-pair
is the vector bundle $\widehat E=(-1)^i \mathrm{Ind}(D)\to \widehat C$ .
\end{defin}

\begin{remark} From a more formal point of view, the study of the family of
Dirac operators $D$ can be approached via the techniques developed by Bismut
in his proofs of the Atiyah-Singer index Theorem for families \cite{Bi}. In
that framework one has to consider the fibration $\pi_{\widehat
C}\colon C\times\widehat C\to \widehat C$ as a family of spin$^c$ manifolds, whose
fibres are precisely $C_{\xi}$. The vector bundle of relative spinors is
identified with $\pi_C^*S$ and we can consider then the coupled relative Dirac
operator $$D\colon \Omega^0(\pi_C^*(S^+\otimes E)\otimes\mathcal{P})\to
\Omega^0(\pi_C^*(S^-\otimes E)\otimes\mathcal{P})\,,$$ whose restriction to
$C_{\xi}$ is  $D_{\xi}$.
\end{remark}

The Nahm transform from $C$ to $\widehat C$ is a procedure which transforms
Hermitian vector bundles with unitary connections on $C$ into Hermitian vector
bundles with unitary connections on $\widehat C$. The main idea relies on the
fact that the index (minus the index) of the family $D$ is a finite
rank vector bundle whenever $\Coker D=0$
($\Ker D=0$). In certain cases this is a consequence of a vanishing
Theorem of Bochner type. Before doing so we introduce some more notation and
recall the Weitzenb\"ock formula.

We recall that the Poincar\'{e} line bundle
$\mathcal{P}\to C\times\widehat C$ is a holomorphic Hermitian line bundle and
that the unitary connection $\nabla_{\mathcal{P}}$ is compatible with the
holomorphic structure.
It is also known that a Hermitian vector bundle $E\to C$ with a
unitary connection $\nabla$ is naturally endowed with a
holomorphic structure since $F^\nabla$ is of type $(1,1)$ (see
\cite[2.1.53]{DK}). Moreover, the spin$^{c}$ Dirac operator
$D_{\xi}$ coincides with the Dolbeault-Dirac operator of
$E\otimes\mathcal{P}_{\xi}$  $$D_{\xi}=\sqrt
2(\bar\partial^*_{E\otimes\mathcal{P}_{\xi}}+\bar\partial_{E\otimes\mathcal{P}_{\xi}}),$$
where $\bar\partial_{E\otimes\mathcal{P}_{\xi}}$ is the
Cauchy-Riemann operator of $E\otimes\mathcal{P}_{\xi}$. Since $C$
is a one dimensional complex manifold the Dolbeault-Dirac operator
$D_{\xi}$ is reduced to

$$D_{\xi}=\sqrt{2}\bar\partial_{E\otimes\mathcal{P}_{\xi}}\colon
\Omega^0(C_{\xi},E\otimes\mathcal{P}_{\xi})\to\Omega^{0,1}(C_{\xi},
E\otimes\mathcal{P}_{\xi})\,.$$

As a consequence of the K\"ahler identities (see \cite{DK}), the Weitzenb\"ock
formula for the Dirac operator $D_\xi$ can be expressed as

\begin{equation}\label{weitzen-el}
D_{\xi}^*D_{\xi}^{}=2\bar\partial_{E\otimes\mathcal{P}_{\xi}}
^*\bar\partial_{E\otimes\mathcal{P}_{\xi}}^{}=\nabla_{\xi}^*\nabla_{\xi}^{}-i\Lambda
F^\nabla\otimes \Id_{\mathcal{P}_{\xi}},
\end{equation} where $i\Lambda F^\nabla$ is the Hermitian endomorphism
of $E$ obtained by contracting $iF^\nabla$ with the K\"ahler form. We have the
following vanishing Theorem.

\begin{thm}\label{vanishing2} Let $(E,\nabla)$ be a pair formed by a Hermitian
vector bundle over $C$ and a unitary connection.
\begin{enumerate}
\item[(i)] If $i\Lambda F^\nabla$ is  non-negative and there exists $x\in C$ such
that $i\Lambda F^\nabla(x)>0$ then $(E,\nabla)$ is an $\IT_0$-pair.
\item[(ii)] If $i\Lambda F^\nabla$ is a non-positive and there exists $x\in C$ such
that $i\Lambda F^\nabla(x)<0$ then $(E,\nabla)$ is an $\IT_1$-pair.
\end{enumerate}

\end{thm}

\begin{proof}
Let us suppose that $i\Lambda F^\nabla<0$. If we apply the Weitzenbock formula
(\ref{weitzen-el}) to a section $s\in\Gamma(C,E\otimes \mathcal{P}_{\xi})$ and
we integrate over $C$ we obtain

\begin{equation}\label{clave}
\Vert D^{}_{\xi}s\Vert^2=\Vert \nabla_{\xi}^{}s\Vert^2-\int_C\langle i\Lambda
F^\nabla s,s\rangle\,\omega \geq 0,
\end{equation}
where $\omega$ is the Riemannian volume element of $C$.
{F}rom relation (\ref{clave}) we obtain $$ D_{\xi}s=0 \Longleftrightarrow
\begin{cases}
(a)\ \nabla_{\xi}^{}s=0\\
(b)\ \langle i\Lambda F^\nabla s,s\rangle=0.
\end{cases}$$
By $(a)$ one sees that  $\langle
s,s\rangle$ is constant; therefore if there exists $x\in C$ such that
$i\Lambda F^\nabla(x)<0$, then $(b)$ implies that $s(x)=0$ and since $\langle
s,s\rangle$ is constant, one has $s=0$ and $(ii)$ is proved.
By Serre duality we have $H^1(C,E)\simeq H^0(C,E^\vee)^*$, and hence the first
statement follows from the second one.
\end{proof}

We can endow the transformed vector bundle of an $\IT$-pair with a
Hermitian metric and a unitary connection in a natural way. This
follows from a rather straightforward application of the theory
for families. We briefly recall the main facts of this
construction following the approach of \cite[Chapter 3]{DK} and
\cite{Bi}.

Let $H^\infty_\pm$ be the space  of $C^\infty$ sections of the vector
bundle $\pi_C^*(S^\pm\otimes E)\otimes\mathcal{P}$ over $C\times \widehat C$.
We may regard $H^\infty_\pm$ as the space  of $C^\infty$ sections over
$\widehat C$ of the infinite dimensional fibre bundles
$\mathcal{H}_\pm^\infty$.
The fibres $\mathcal{H}_{\pm,\,\xi}^\infty$ are the sets of $C^\infty$
sections over $C_{\xi}$ of $\pi_C^*(S^\pm\otimes E)\otimes\mathcal{P}$. Since
$\pi_C^*(S^\pm\otimes E)\otimes\mathcal{P}$ is a Hermitian vector bundle, and
the fibres $C_{\xi}$ of the projection $\pi_{\widehat C}\colon C\times\widehat
C\to\widehat C$ carry a natural volume element $\omega$; we can define the Hermitian metric
\begin{equation}\label{met-herm}
\langle h_1,h_2\rangle_{\pi_{\widehat C}}=\int_{C_{\xi}}\langle
h_1,h_2\rangle\,\omega\,,
\end{equation} on $\mathcal{H}_{\pm,\,\xi}^\infty$,
 We then have the Hilbert bundles $\mathcal{H}_{\pm}$ whose
fibres $\mathcal{H}_{\pm,\,\xi}$ are the $L^2$-completion of
$\mathcal{H}_{\pm,\,\xi}^\infty$ with respect to this metric.

Let $\nabla^1$ be the connection on $\pi_C^*(S^\pm\otimes
E)\otimes\mathcal{P}$ obtained from $\nabla_S$, $\nabla$ and
$\nabla_{\mathcal{P}}$. Now we define a connection $\widetilde\nabla$ on
$\mathcal{H}_\pm^\infty$ as follows $$\widetilde\nabla_D h=\nabla_{D^H}^1h,\ \
\ \  \text{for every} \ \  D\in\mathfrak{X}(\widehat C),\ h\in H^\infty_\pm, $$
where $D^H$ is  the natural lift of the vector field $D$ from $\widehat C$ to
$C\times\widehat C$. It is easy to check that $\widetilde\nabla$ is
a flat connection.

If $(E,\nabla)$ is an $\IT_i$-pair, then the regularity Theorem for elliptic
operators implies that $\widehat E$ is, according to the parity of the index
$i$, a subbundle of $\mathcal{H}_\pm^\infty$, and hence there is a naturally
induced metric on $\widehat E$. We also have a natural unitary connection
$\widehat\nabla$ induced by the ambient connection $\widetilde\nabla$ and the
orthogonal projection $P$ onto $\widehat E$, that is
$$\widehat\nabla=P\circ\widetilde\nabla.$$ Let us recall that Hodge theory
provides an explicit formula for the projector $P$. Indeed, if $(E,\nabla)$ is
$\IT_0$ then for every $\xi\in\widehat C$ we have
$$P_{\xi}=\Id-D^*_{\xi}G_{\xi}^{}D_{\xi}^{},$$ where $G_{\xi}^{}$ is the Green
operator of $D_{\xi}^{}D^*_{\xi}$. A similar formula holds in the case of an
$\IT_1$ pair.

\begin{defin}
Let $(E,\nabla)$ be an $\IT$-pair. The pair $(\widehat E,\widehat\nabla)$
is called the Nahm transform of $(E,\nabla)$ and is denoted by $\mathcal{N}(E,\nabla)$.
\end{defin}

\begin{rem}\label{gauge} If $\nabla$ and $\nabla^\prime$ are gauge
equivalent unitary connections, it follows from the very definition of the Nahm
transform that  $\widehat\nabla$ and $\widehat{\nabla^\prime}$ are also gauge
equivalent unitary connections.
\end{rem}

The following is an easy consequence of the flatness of $\widetilde\nabla$.

\begin{prop}\label{curvat}
Let $(\widehat E,\widehat\nabla)$ be the Nahm transform of an $\IT_i$ pair.
The curvature of $\widehat\nabla$ is given by $$F^{\widehat\nabla}=
P\circ(\widetilde\nabla P\wedge \widetilde\nabla P)\circ P.$$ Moreover, we can
express the curvature in terms of the Green operator as follows
\begin{align*}
F^{\widehat\nabla}= P\circ(\widetilde\nabla D^*\circ G\wedge \widetilde\nabla
D)\circ P,\ \ \ \ \ \mathrm{if}\ \ E\ \mathrm{is}\ \IT_0.
\end{align*} A similar expression holds in the case of an $\IT_1$ pair.
\qed\end{prop}

We study now the Nahm transform of a connection with constant central
curvature.   Since all the line bundles $\Pc_\xi$ are flat they are trivial as
smooth bundles and we may consider the connection $\nabla_\mathcal{P}$ of the
Poincar\'{e} line bundle as a family of connections $\overline\nabla_{\xi}$ on the
trivial line bundle. In the same way if $E\to C$ is a Hermitian vector bundle
with a unitary connection $\nabla$ then we have a family of connections
$\nabla_\xi$ on $E$ and the family of Dolbeault-Dirac operators $D_\xi$
considered above act now in the same vector bundle
$$\overline\partial^{\nabla_\xi}\colon \Omega^0(E)\to\Omega^{0,1}(E).$$ Using
a flat holomorphic coordinate $z$ on $C$ and the flat coordinate $w$ which it
induces on $\widehat C$ we have
$$\overline\partial^{\nabla_\xi}=\overline\partial^{\nabla}+\pi\,wd\bar
z\otimes \Id_{E},$$ which clearly shows that this family depends
holomorphically on $w\in\widehat C$.

The triviality of the holomorphic tangent bundle of $C$ allows to identify
$\Omega^{0,1}(C)$ with $\Omega^0(C)$ by contraction with a global
anti-holomorphic vector field $\overline V$. Since the metric on $C$ is flat,
we can choose $\overline V$ such that it is a parallel vector field whose
pointwise norm is equal to 1. We define  the operator $$
\mathcal{D}_\xi=i_{\overline
V}\dbar^{\nabla_\xi}\colon\Omega^0(E)\lra\Omega^0(E)\,. $$

\begin{lemma}
\label{commutator} The curvature of $\nabla_\xi$ is related to the
operator $\mathcal{D}_\xi$ by the formula $$ i\Lambda
F^{\nabla}=i\Lambda
F^{\nabla_\xi}=2[\mathcal{D}_\xi,\mathcal{D}_\xi^\ast].$$
\end{lemma}

Let us recall that the flat metric of $C$ induces in a natural way a flat
metric on $\widehat C$.
In the following Theorem we consider unitary
connections of constant central curvature on $C$ and $\widehat C$ with respect
to these metrics.

\begin{thm}\label{eh}
Let $\nabla$ be a connection on $E$ with constant central curvature
with factor $\lambda\in\mathbb{R}$, that
is $i\Lambda F^\nabla=\lambda\, \mathrm{Id}_E$,
where $\lambda=2\pi\mu(E)$ and $\mu(E)$ is the slope of $E$.

\begin{enumerate}
\item If $\mathrm{deg}(E)>0$ then $(E,\nabla)$ is an $\IT_0$ pair and $\widehat\nabla$
is a connection on $\widehat E$  with constant central curvature with factor
$\hat\lambda=-\frac{2\pi}{\mu(E)}$.

\item If $\mathrm{deg}(E)<0$ then $(E,\nabla)$ is an $\IT_1$ pair and $\widehat\nabla$
is a connection on $\widehat E$ with constant central curvature with factor
$\hat\lambda=-\frac{2\pi}{\mu(E)}$.
\end{enumerate}

\end{thm}

\begin{proof} We shall only prove the first case since the second one can be
dealt with in a similar way.

It is well-known (\cite{Do} and \cite{NS}) that since  $\nabla$ has constant
central curvature $E$ must be polystable. The condition $\mathrm{deg}(E)>0$
implies, due to Proposition \ref{presstab}, that $(E,\nabla)$ is an $\IT_0$
pair. All the operators $\mathcal{D}_\xi$ act on $\Omega^0(E)$, therefore the
bundle of kernels $\widehat E$ is a finite rank subbundle of the trivial
Hilbert bundle $\mathcal{H}_+\to\widehat C$ introduced above and
$P\colon\mathcal{H}_+\to\widehat E$ is the orthogonal projection. Then we have
$$\widehat\nabla=P\circ\widetilde\nabla,$$ where $\widetilde\nabla$ is the
natural flat connection on $\mathcal{H}_+$. Taking into account the above
identifications, the curvature of the connection $\widehat\nabla$ of $\widehat
E$, given in Proposition \ref{curvat}, can be expressed as
\begin{equation}\label{curv-f}
F^{\widehat\nabla}=P_\xi\circ(\widetilde\nabla\mathcal{D}^*_\xi\circ
G_\xi\wedge\widetilde\nabla\mathcal{D}_\xi)\circ P_\xi,
\end{equation} where
$G_\xi$ is the Green operator of
$\mathcal{D}_\xi\mathcal{D}_\xi^*$.

As we mentioned above, we can choose a flat holomorphic coordinate $z$ on $C$
such that the K\"ahler form is expressed as
$$\omega=\frac{i}{2}dz\wedge d\bar
z\,.$$

 Therefore, locally we may take $\overline V=\frac{\partial}{\partial\bar
z}$. This implies that $$\mathcal{D}_\xi=\mathcal{D}_0+\pi w \Id_E.$$

It is clear now that $\widetilde\nabla\mathcal{D}_\xi=\pi dw\otimes \Id_E$ and
$\widetilde\nabla\mathcal{D}_\xi^*=\pi d\bar w\otimes \Id_E$ which upon
substitution in (\ref{curv-f}) gives $$ F^{\widehat\nabla}=\pi^2\,P_\xi\circ
G_\xi\circ P_\xi\, d\bar{w}\wedge dw, $$ where we have used the fact that the
identity operator commutes with the Green's operator $G_\xi$. We then have to
prove that for every $u\in\ker\mathcal{D}_\xi$ one has $$ G_\xi u=\alpha
u+v\,, $$ where $\alpha$ is a constant and $v\in (\ker
{\mathcal{D}}_\xi)^\perp$. To see this suppose that $$ G_\xi u=u'+v
\;\;\;\;\mbox{for}\;\; u'\in\ker {\mathcal{D}}_\xi\;\;\;\mbox{and}\;\;\; v\in
(\ker {\mathcal{D}}_\xi)^\perp. $$ Operating by
$G_\xi^{-1}=\mathcal{D}_\xi\mathcal{D}_\xi^\ast$ we obtain
\begin{equation}
u=\mathcal{D}_\xi\mathcal{D}_\xi^\ast
u'+\mathcal{D}_\xi\mathcal{D}_\xi^\ast v. \label{g-1}
\end{equation}
But by hypothesis $[\mathcal{D}_0,\mathcal{D}_0^\ast]=\frac{\lambda}{2}\,
\Id_E$, therefore
$[{\mathcal{D}}_\xi,{\mathcal{D}}_\xi^\ast]=\frac{\lambda}{2}\,\Id_E$ and,
since $\mathcal{D}_\xi u^\prime=0$, equation (\ref{g-1}) becomes $$
u-\frac{\lambda}{2}\, u'=\mathcal{D}_\xi\mathcal{D}_\xi^\ast v. $$ Now
$\mathcal{D}_\xi\mathcal{D}_\xi^\ast v\in (\ker {\mathcal{D}}_\xi)^\perp$,
since for every $u_1\in\ker {\mathcal{D}}_\xi$,
\begin{eqnarray}
(\mathcal{D}_\xi\mathcal{D}_\xi^\ast v, u_1)&=&(\frac{\lambda}{2}\, v+
{\mathcal{D}}_\xi^\ast {\mathcal{D}}_\xi v,u_1)\nonumber\\
                  &=&(\frac{\lambda}{2}\, v,u_1)+({\mathcal{D}}_\xi v,
                  {\mathcal{D}}_\xi u_1)\nonumber\\
                  &=&0. \nonumber
\end{eqnarray}
Thus $u-\frac{\lambda}{2}\, u'\in \ker {\mathcal{D}}_\xi\cap (\ker
{\mathcal{D}}_\xi)^\perp=\{0\}$. Hence $u'=2\lambda ^{-1} u$, concluding that
$$ F^{\widehat\nabla}=-2\pi^2\lambda^{-1} dw\wedge
d\bar{w}=-\frac{(2\pi)^2}{i\lambda} \frac{i}{2}dw\wedge
d\bar{w}=-\frac{(2\pi)^2}{i\lambda}\widehat\omega\otimes \Id_{\widehat E},$$
where $\widehat\omega$ is the K\"ahler form of $\widehat C$. Therefore
$i\Lambda F^{\widehat\nabla}=-\frac{2\pi}{\mu(E)}\, \Id_{\widehat E}$ as
required.
\end{proof}

\subsection{Compatibility between the Fourier-Mukai and
Nahm transforms, functoriality and invertibility}\label{comp}

Let $E\to C$ be a Hermitian vector bundle endowed with a unitary connection
$\nabla$. As we have seen  the spin$^c$ Dirac operator $D_{\xi}$ is
identified with the Dolbeault-Dirac operator of
$E\otimes\mathcal{P}_{\xi}$. Hodge theory and the Dolbeault isomorphism give that
\begin{align}\label{dol1} \Ker D_{\xi}&\simeq H^{0}(C_{\xi},E\otimes\mathcal{P}_{\xi})\\
\label{dol2} \Coker D_{\xi} &\simeq
H^{1}(C_{\xi},E\otimes\mathcal{P}_{\xi}).
\end{align}

If we suppose that $E$ is $\IT_i$ with respect to the
Fourier-Mukai transform $\mathcal{S}$, then the isomorphisms
(\ref{dol1}) and (\ref{dol2}) mean that $(E,\nabla)$ is an
$\IT_i$-pair with respect to the Nahm transform. By \cite[Theorem
2]{BBH1} or \cite[Theorem 3.2.8]{DK} we have a natural $C^\infty$
vector bundle isomorphism induced by Hodge theory
$$\phi_E\colon\widehat E\iso\mathcal S^i (E).$$ Moreover, we have the
following.

\begin{thm}\label{compatibility} Let $E_1$, $E_2$ be Hermitian vector
bundles over $C$ endowed with unitary connections $\nabla_1$,
$\nabla_2$ such that $(E_1,\nabla_1)$, $(E_2,\nabla_2)$ are
$\IT_i$-pairs with respect to the Nahm transform. Then we have
\begin{enumerate}\item The connections $\widehat\nabla_1$,
$\widehat\nabla_2$ are compatible with the holomorphic structures
of $\mathcal S^i (E_1)$, $\mathcal S^i (E_2)$, respectively.
\item For every holomorphic morphism $\Phi\colon E_2\to E_1$
we have an induced
holomorphic
morphism
$\mathcal{N}(\Phi)\colon\widehat E_2\to \widehat E_1$ and a
commutative diagram
\begin{equation}
\xymatrix{{\widehat E_2}\ar[r]^{\phi_{E_2}}\ar[d]_{\mathcal{N}(\Phi)} &
{\mathcal{S}(E_2)}\ar[d]^{\mathcal{S}(\Phi)}\\
{\widehat E_1}\ar[r]^{\phi_{E_1}} & {\mathcal{S}(E_1)} }
\label{comm}
\end{equation}
\end{enumerate}
\end{thm}
\begin{proof}
The Poincar\'{e} bundle $\mathcal{P}\to C\times\widehat C$ is a
holomorphic line bundle and the connection $\nabla_{\mathcal{P}}$
is compatible with the holomorphic structure. This implies that
the families of Dirac operators $D_{\xi}^{}$, $D_{\xi}^{*}$ vary
holomorphically with $\xi\in\widehat C$. The first statement
follows now by a standard argument concerning holomorphic
families, see \cite[Theorem 3.2.8]{DK}.

Since $\Phi$ is holomorphic, the second statement follows
immediately in the $\IT_0$ case, because the fibers of the Nahm
transforms are given by the kernels of the Cauchy-Riemann
operators. In the $\IT_1$ case the fibers $\widehat E_{2,\xi}$,
$\widehat E_{1,\xi}$ of the Nahm transform at $\xi\in\widehat C$
are given by the cokernels of the Dirac operators $D_{2,\xi}$,
$D_{1,\xi}$ and we have
\begin{align*}\Coker D_{2,\xi}&=\Ker D_{2,\xi}^*= \Ker
\overline\partial^{\nabla_{2,\xi}^*}\\
\Coker D_{1,\xi}&=\Ker D_{1,\xi}^*= \Ker
\overline\partial^{\nabla_{1,\xi}^*}
\end{align*} Now $\Phi$ induces a morphism from $\Ker
\overline\partial^{\nabla_{2,\xi}^*}$ to
$\Omega^{0,1}(C,E_1\otimes\mathcal{P}_\xi)$ and composing it with
the orthogonal projection onto $\Ker
\overline\partial^{\nabla_{1,\xi}^*}$ we get a morphism
$$\mathcal{N}(\Phi)_\xi\colon \Ker
\overline\partial^{\nabla_{2,\xi}^*}\to  \Ker
\overline\partial^{\nabla_{1,\xi}^*}$$  which by Hodge theory is
the unique one that renders commutative the following diagram
$$\xymatrix{{\widehat
E_{2,\xi}}\ar[r]^(.4){\phi_{E_2,\xi}}\ar[d]_{\mathcal{N}(\Phi)_\xi}
&
{ \mathcal{S}(E_2)_\xi} \ar[d]^{\mathcal{S}(\Phi)_\xi}&{\hskip-1cm=H^1(C,E_{2,\xi})}\\
{\widehat E_{1,\xi}}\ar[r]^(.4){\phi_{E_{1,\xi}}} & {
\mathcal{S}(E_1)_{\xi}}&{\hskip-1cm=H^1(C,E_{1,\xi})} }$$ Since
$\mathcal{S}(\Phi)$ is a vector bundle morphism and $\phi_{E_2}$, $\phi_{E_1}$
are $C^\infty$ vector bundle isomorphisms, we conclude that
$\mathcal{N}(\Phi)$ is also a $C^\infty$ vector bundle morphism and we have
the commutative diagram \eqref{comm}. Moreover, $\mathcal{S}(\Phi)$ is a
holomorphic morphism and by the first part of the Theorem we have the
compatibility between the connections $\widehat\nabla_2$, $\widehat\nabla_1$
and the holomorphic structures of $\mathcal{S}(E_2)$, $\mathcal{S}(E_1)$,
respectively. These facts imply that $\mathcal{N}(\Phi)$ is a holomorphic
morphism.

\end{proof}

If $h$ is an Hermitian metric on a $C^\infty$ vector bundle $E$ then
$\mathcal{A}(E,h)$ will denote the space of unitary connections which are
compatible with $h$. On the other hand we will denote by $\mathcal{C}(E)$ the
set of holomorphic structures on $E$. It is well known that there is an
identification $$\xymatrix{{\mathcal{A}(E,h)}\ar[r]^{\ \ \sim} &
{\mathcal{C}(E)}}$$ which associates to $(E,\nabla)$ the holomorphic vector
bundle $\mathcal{E}=(E,\bar\partial^\nabla)$, the inverse
correspondence being given by the map which associates to every
holomorphic bundle $\mathcal{E}=(E,\bar\partial^\nabla)$ the unique
connection $\nabla$ compatible with the complex structure and the
Hermitian metric. We can rephrase the preceding
Theorem by saying that the Nahm transform and the Fourier-Mukai transform are
compatible with this identification. That is to say, the following diagram is
commutative $$\xymatrix{{\mathcal{A}(E,h)}\ar[r]^{\ \
\sim}\ar[d]_{\mathcal{N}} &
{\mathcal{C}(E)}\ar[d]^{\mathcal{S}}\\
{\mathcal{A}(E,h)}\ar[r]^{\ \ \sim} & {\mathcal{C}(E).} }$$

The curve $C$ and its dual elliptic curve $\widehat C$ are in a
symmetrical dual relation with one another (see \cite[Section
3.3.2] {DK}). That is, $C$ parametrizes the flat Hermitian line
bundles over $\widehat C$, therefore $\widehat{\widehat C}\simeq
C$. Moreover, the restriction of the dual of the Poincar\'{e} line
bundle $\mathcal{P}^\vee$ to the slice $\widehat C_{x}$ endowed
with the restriction of the connection $\nabla_{\mathcal{P^\vee}}$
is isomorphic, as a Hermitian line bundle with connection, to the
flat Hermitian bundle corresponding to $x$. We can hence apply the
Nahm construction in order to transform Hermitian vector bundles
with connection over $\widehat C$ into Hermitian vector bundles
with connection over $C$.

Let  $\nabla$ be a connection with constant central curvature different from
zero on a bundle  $E$ over $C$, and let $\mathcal{E}=(E,\bar\partial^\nabla)$
be the corresponding holomorphic vector bundle; then $\deg(\mathcal{E})\neq
0$. The isomorphisms (\ref{dol1}) and (\ref{dol2}) imply that $\mathcal{E}$ is
$\IT_i$ with respect to the Fourier-Mukai transform. Let
$\widehat{\mathcal{E}}=\mathcal S^i  (E)$ be  its unique transform. It is well
known, see \cite{Mu}, that $\widehat{\mathcal{E}}$ is $\IT_{1-i}$ and that
there is an isomorphism of holomorphic vector bundles
\begin{equation}\label{iso}
\widehat{\widehat{\mathcal E}}= {\widehat{\mathcal{S}}}^{\,1-i}(
\mathcal S^i(E))\simeq \mathcal{E}.
\end{equation}

By Theorem \ref{eh} $\widehat\nabla$ is a constant central curvature
connection on $\widehat E$, and hence we can apply to it the Nahm transform to
obtain $(\widehat{\widehat{E}},\widehat{\widehat{\nabla}})$. By (\ref{iso}) we
have an isomorphism $$\widehat{\widehat{E}}\simeq E\,.$$

Moreover,  Theorem \ref{compatibility} implies that
$\widehat{\widehat{\nabla}}$ is compatible with the holomorphic
structure of $\mathcal{E}$, and therefore  by
the results of Donaldson \cite{Do}, which in particular extend the
theorem of Narasimhan and Seshadri \cite{NS} to genus one,
we have the following.

\begin{thm} If $\nabla$ is a connection  with constant central
curvature different from zero on $E$ then $\widehat\nabla$ is a
connection with constant central curvature on the bundle $\widehat
E$, and there is a natural isomorphism
$$(\widehat{\widehat
E},\widehat{\widehat{\nabla}})\simeq(E,\nabla).$$
\end{thm}

Let $\mathcal{A}_c(E,h)\subset \mathcal{A}(E,h)$ be the subspace
of constant central curvature connections  and
let $\mathcal{C}_{ps}(E)\subset\mathcal{C}(E)$  be the subspace of polystable
holomorphic structures on the $C^\infty$ bundle $E$.
We have the Donaldson--Narasimhan--Seshadri  correspondence (the curve version of the
Hitchin--Kobayashi correspondence)
$$
{\mathcal{A}_c(E,h)}\xrightarrow{\ D\ } {\mathcal{C}_{ps}(E)}\,.
$$

The content of the preceding Theorem can be summarized by saying that the Nahm
transform and the Fourier-Mukai transform are compatible with
the Donaldson--Narasimhan--Seshadri
correspondence. That is, the following diagram is
commutative $$\xymatrix{{{\mathcal{A}}_c(E,h)}\ar[r]^{D}\ar[d]_{\mathcal{N}} &
{{\mathcal{C}}_{ps}(E)}\ar[d]^{\mathcal{S}}\\
{{\mathcal{A}}_c(\widehat E,\hat h)}\ar[r]^{D} &
{{\mathcal{C}}_{ps}(\widehat E)\,.} }$$

These correspondences descend to the quotients by the corresponding gauge
groups, giving a commutative diagram of correspondences between the associated
moduli spaces. First, the Donaldson--Narasimhan--Seshadri correspondence is
well known to descend to moduli spaces, see \cite[Chapter VII]{kobayashi}. The
descent  for the Nahm transform follows from Remark \ref{gauge} and  for the
Fourier-Mukai transform is a consequence of its functoriality.

\section{Fourier-Mukai transforms for holomorphic triples}\label{transtrip}

\subsection{Holomorphic triples}
A holomorphic triple over a smooth connected curve $C$ is by definition a
triple $T= (E_1, E_2, \Phi)$ where $E_i$, $i=1,2$ are holomorphic vector
bundles and $\Phi \in \Hom_C (E_2,E_1)$. Let
$n_i$ and  $d_i$ be the rank and degree of $E_i$ for $i=1,2$. We say that
the triple $T$ is of type $(n_1,n_2,d_1,d_2)$. There is a notion of stability
for triples which depends on a real parameter $\alpha $ (see \cite{BG} for
details). The $\alpha $-degree of $T$ is defined by $$\deg_\alpha (T) = \deg
(E_1\oplus E_2) +n_2 \alpha $$ \noindent and the $\alpha$-slope is by
definition $$ \mu_\alpha (T) = \frac{\deg_\alpha (T)}{n_1+n_2}.$$ The
stability condition is defined in a similar way as the slope stability for
vector bundles, precisely: $T=(E_1, E_2, \Phi)$ is $\alpha$-stable (resp.
$\alpha$-semistable) if for every non-trivial subtriple $ T^\prime =
(E_1^\prime, E_2^\prime, \Phi^\prime )$ we have
\begin{equation*}
 \mu_\alpha (T^\prime) <  \mu_\alpha (T) \qquad  (\text{resp.} \leq ).
\end{equation*}
Here a subtriple means a triple $T^\prime = (E_1^\prime, E_2^\prime,
\Phi^\prime   )$ and injective homomorphisms $\gamma_1, \gamma_2$ of sheaves
such that the following diagram commutes
\begin{equation*}
\xymatrix{ E_2^\prime  \ar[r]^{\Phi^\prime} \ar[d]_{\gamma_2} &
E_1^\prime \ar[d]^{\gamma_1} \\
E_2  \ar[r]^\Phi & E_1. }
\end{equation*}

Most of the properties which are valid for stable bundles carry along to stable
triples. We denote the moduli space of S-equivalence classes of
$\alpha$-semistable triples of type $(n_1,n_2,\break d_1,d_2)$ by $\mathcal N_\alpha
(n_1,n_2,d_1,d_2)$ or simply by $\mathcal N_\alpha $ if there is no need to
specify the topological invariants. $\mathcal N_\alpha^s (n_1,n_2,d_1,d_2)$
denotes the moduli space of $\alpha$-stable triples.

An important feature is that the stability condition gives bounds on the
range of the parameter $\alpha $. More precisely, if $n_1 \neq n_2$
and $T=(E_1, E_2, \Phi)$ is $\alpha$-stable of type $(n_1,n_2,d_1,d_2)$ then
necessarily $$ 0\leq\alpha_m \leq \alpha \leq \alpha_M $$ where $\alpha_m =
\mu_1 - \mu_2$ and $\alpha_M = \left( 1 +\frac{n_1+n_2} {|n_1 - n_2|}\right)
(\mu_1 - \mu_2)$ (see \cite{BG}). In the case $n_1=n_2$, $\alpha $ ranges in
$[\alpha_m, \infty)$; we will write in this case, $\alpha_M =\infty$.
The interval $(\alpha_m, \alpha_M)$ is divided into a finite number of
subintervals  determined by values of the parameter
for which strict  semistability may occur. The stability criteria for two
values of $\alpha$ lying between two consecutive critical values are
equivalent (and therefore the corresponding moduli spaces are isomorphic). As
in \cite{BGG1} we shall denote by $\alpha_L$ the largest critical
value, in particular when $ \alpha_L< \alpha < \alpha_M $ all the
moduli spaces $\mathcal N_\alpha$ are isomorphic.

We use freely the terminology and results of \cite{BG}.  Corollary 3.6,
Proposition 3.17, Corollaries 3.19 and 3.20 and Lemma 4.6 of \cite{BG} are
particularly useful for the  understanding of this paper.

Now we recall how  holomorphic triples on an elliptic curve $C$ are
related to $SU(2)$-equivariant bundles on the elliptic surface $C\times
\mathbb P^1$. In what follows we shall only deal with $SU(2)$-equivariant
bundles $E$ which admit a $C^\infty$ $SU(2)$-equivariant decomposition
of the  type

\begin{equation}\label{type}
E = p^* E_1 \oplus (p^* E_2 \otimes q^* H^{\otimes 2}),
\end{equation}
where $p$, $q$ are the canonical projections of $C\times \mathbb
P^1$ onto its factors and $ H$ is the $C^\infty$ line bundle over
$\mathbb{P}^1$ with first Chern number equal to $1$.

In the following, if not otherwise stated, an $SU(2)$-equivariant bundle will
always mean an holomorphic bundle over $C\times \mathbb P^1$,
$SU(2)$-equivariant, of type given in (\ref{type}).

We shall need the following formulation of Proposition 2.3 in
\cite{BG}.

\begin{prop}\label{invv} Let $C$ be a smooth connected curve, then
\begin{itemize}
\item[(i)] There is a one-to-one
correspondence between $SU(2)$-equivariant holomorphic vector bundles $E$ of
type (\ref{type}) and holomorphic extensions over $C\times \mathbb P^1$ of the
form $$ 0\to p^* E_1 \to E \to p^* E_2 \otimes q^* \mathcal O_{\mathbb{P}^1}
(2)\to 0$$ \noindent where $E_1$, $E_2$ are holomorphic vector bundles on $C$.
Here $\mathcal O_{\mathbb{P}^1}(2)$ is the unique line bundle of degree $2$
over $\mathbb P^1$.
\item[(ii)] There is a (non-unique) functorial correspondence between such
extensions and elements of $\Hom_C (E_2,E_1)$ and it is given by a functorial
isomorphism $$\sigma_C\colon \Ext^1_{C\times \mathbb P^1} (p^*E_2 \otimes
q^*\mathcal O_{\mathbb P^1} (2), p^*E_1) \simeq \Hom_C (E_2, E_1)$$ induced by
the choice of a trace isomorphism $\tr\colon H^1(\mathbb P^1, \mathcal
O_{\mathbb P^1} (-2))\iso \C$.
\end{itemize}
\qed\end{prop}

\begin{proof}
A proof of $(i)$ and $(ii)$ can be found in \cite[Proposition 3.9]{G1} and
\cite[Proposition 2.3]{BG}.
We recall that there is a natural isomorphism (see for instance \cite{H})
\begin{equation}
\Ext^1_{C\times \mathbb P^1} (p^*E_2 \otimes q^*\mathcal
O_{\mathbb P^1} (2), p^*E_1)\iso \Hom_{D(C\times \mathbb P^1)}(p^*E_2 \otimes q^*\mathcal
O_{\mathbb P^1} (2), p^*E_1[1])\,.
\label{ext1}
\end{equation}
Now we have
\begin{align*}\Hom_{D(C\times \mathbb P^1)}(p^*E_2 \otimes q^*\mathcal O_{\mathbb P^1} (2),
p^*E_1[1])& \iso  \Hom_{D(C\times \mathbb P^1)} (p^*E_2 , p^*E_1\otimes
q^*\mathcal O_{\mathbb P^1} (-2) [1]) \\
&\iso \Hom_{D(C)}(E_2 , \Rb p_*(p^*E_1\otimes q^*\mathcal
O_{\mathbb P^1} (-2)) [1]) \\
& \iso \Hom_{D(C)}(E_2 , E_1\otimes \Rb p_*(q^*\mathcal
O_{\mathbb P^1} (-2)) [1]) \\
&\iso \Hom_{D(C)}(E_2 , E_1\otimes_{\C} \Rb \Gamma(\mathbb P^1,\mathcal
O_{\mathbb P^1} (-2)) [1])\\
& \iso  \Hom_{D(C)}(E_2 , E_1\otimes_{\C} H^1(\mathbb P^1,\mathcal O_{\mathbb
P^1} (-2)))
\end{align*} where the second isomorphism is adjunction between direct and
inverse images, the third is the projection formula, the fourth is base-change
in the derived category and the last is due to the fact that since $H^0(\mathbb
P^1,\mathcal O_{\mathbb P^1} (-2))=0$, then $\Rb \Gamma(\mathbb P^1,\mathcal
O_{\mathbb P^1} (-2))\iso H^1(\mathbb P^1,\mathcal O_{\mathbb P^1} (-2)) [-1]$
in the derived category. Composition with a trace map $\tr\colon H^1(\mathbb
P^1, \mathcal O_{\mathbb P^1} (-2))\iso \C$ gives the isomorphism $$
\sigma_C\colon \Ext^1_{C\times \mathbb P^1} (p^*E_2 \otimes q^*\mathcal
O_{\mathbb P^1} (2), p^*E_1) \simeq \Hom_C (E_2, E_1) $$ of the statement.

\end{proof}

\begin{remark} We can describe quite easily in an explicit form the inverse
isomorphism $\sigma_C^{-1}$. The inverse of the trace $\tr^{-1}\colon\C\iso
H^1(\mathbb P^1,\mathcal O_{\mathbb P^1} (-2))$ defines an element of
$H^1(\mathbb P^1,\mathcal O_{\mathbb P^1} (-2))$ and via the isomorphism
\begin{align*}
\Hom_{D(\mathbb P^1)}(\mathcal O_{\mathbb P^1}(2), \mathcal O_{\mathbb
P^1}[1])& \iso\Ext^1_{\mathbb P^1}(\mathcal O_{\mathbb P^1}(2), \mathcal
O_{\mathbb P^1})\iso \Ext^1_{\mathbb P^1}(\mathcal O_{\mathbb P^1}, \mathcal
O_{\mathbb P^1}(-2)) \\
& \iso H^1(\mathbb P^1,\mathcal O_{\mathbb P^1} (-2))
\end{align*}
induces a morphism $\tr^{-1}\colon\mathcal O_{\mathbb P^1}(2)\to \mathcal
O_{\mathbb P^1}[1]$ in the derived category. Thus, given a morphism $\Phi\colon
E_2\to E_1$, one finds that $\sigma_C^{-1}(\Phi)$ is the element of $\Ext^1_{C\times \mathbb
P^1} (p^*E_2 \otimes q^*\mathcal O_{\mathbb P^1} (2), p^*E_1)$ corresponding
to the morphism $$ p^*(\Phi) \otimes q^*(\tr^{-1})\colon p^*E_2 \otimes
q^*\mathcal O_{\mathbb P^1} (2) \to p^*E_1[1] $$ by the isomorphism
\eqref{ext1}. \label{explicitPhi}
\end{remark}

\begin{remark}
Two triples $(E_1,E_2,\Phi)$ and $(E_1,E_2,\lambda\Phi)$ ($\lambda\neq 0$)
define different extensions though  the same holomorphic bundle. However,
they define different $SU(2)$-equivariant holomorphic vector bundles
(see \cite{G}, \cite{BG}), because extensions correspond to $SU(2)$-equivariant
holomorphic vector bundles and not merely to holomorphic vector bundles.
\end{remark}

The correspondence in Proposition  \ref{invv} also preserves
stability. Let $\omega_\alpha $ be the K\"ahler class over
$X\times {\mathbb P}^1$ defined by $\omega_\alpha  =
\frac{\alpha}{2} p^*\omega_C + q^*\omega_{\mathbb P^1}$, with
$\alpha \in {\mathbb R}^+ $. The following result is proved in
\cite{BG}.

\begin{thm}\label{dimred} Let $T=(E_1,E_2,\Phi)$ be an holomorphic
triple over a smooth connected curve $C$ and let $E$ be the
holomorphic $SU(2)$-equivariant bundle $C\times \mathbb P^1$
defined in Proposition \ref{invv}. Then, if $E_1$ and $E_2$ are
not isomorphic, $T$ is $\alpha $-stable if and only if  $E$ is
slope-stable with respect to the K\"ahler form $\omega_\alpha $.
In the case $E_1 \simeq E_2$ then $T$ is $\alpha $-stable if and
only if $\Phi\neq 0$, $E_1\simeq E_2$ is stable and $E$  decomposes as a
direct sum
$$ E\simeq (p^* E_1 \otimes q^* \mathcal O_{\mathbb P^1} (1) )
\oplus ( p^* E_2 \otimes q^* \mathcal O_{\mathbb P^1} (1) )\,.
$$
 \end{thm}

\begin{remark}\label{family}
The proof of statement (ii) in Proposition \ref{invv} that we have just given
above shows that the correspondence between triples and equivariant bundles
also extends to families. Indeed, families of (stable) triples correspond
functorially to families of $SU(2)$-equivariant (stable) bundles.  This
implies that the moduli space $\mathcal N_\alpha $ of $\alpha $-stable triples
(of a given topological type) over an elliptic curve corresponds, via the
canonical isomorphism of Proposition \ref{invv}, to a component  of
the moduli space $\mathcal
M^{SU(2)}_\alpha $ of $SU(2)$-equivariant bundles (defined by the
lift of the $SU(2)$ action determined by (\ref{type}))
stable with respect to the K\"ahler form $\omega_\alpha $.
Therefore  we have a canonical identification $$
 \mathcal N_\alpha \iso  \mathcal M^{SU(2)}_\alpha.
$$
\end{remark}

\subsection{Fourier-Mukai transforms for triples}\label{fm-triple}

We begin by briefly recalling the main properties of relative Fourier-Mukai
transform in the case of a trivial elliptic fibration over the projective line.

The corresponding functor is then
\begin{gather}
\mathcal S_{\mathbb P^1}\colon D(C\times \mathbb P^1) \to D(\widehat C \times
\mathbb P^1)\\
\mathcal S_{\mathbb P^1} (-) = \Rb \pi_{\widehat C\times \mathbb P^1,*}
(\pi_{C\times \mathbb P^1}^* (-) \otimes  \pi_{C\times \widehat C}^*(\Pc)) \notag
\end{gather}
\noindent where $\pi_{C\times {\mathbb P^1}}$,  $\pi_{\widehat C\times
{\mathbb P^1}}$ and $\pi_{C\times \widehat C}$ are the canonical projections of
$C\times \widehat C\times \mathbb P^1$ onto its factors. As in Section \ref{section:FM},
this functor is invertible (see for instance \cite{HP,Mu1}).

We also know that the relative Fourier-Mukai transform is compatible with base-change in
the derived category \cite{HP}. In particular, for vector bundles $E$ in $D(C)$ and $F$
in $D(\mathbb P^1)$ the base change isomorphism can be described as follows. Let us denote
by $\hat p$, $\hat q$ the projections of $\widehat C\times \mathbb P^1$ onto its factors. Then
\begin{equation}
\begin{aligned}
\mathcal S_{\mathbb P^1}(p^* E\otimes q^* F)&\iso  \Rb \pi_{\widehat C\times \mathbb P^1,*}
(\pi_{C\times \mathbb P^1}^* (p^\ast E\otimes q^* F) \otimes  \pi_{C\times \widehat C}^*(\Pc))\\
& \iso \Rb \pi_{\widehat C\times \mathbb P^1,*}
(\pi_{\widehat C\times \mathbb P^1}^*(\hat q^* F) \otimes
\pi_{C\times \widehat C}^*(\pi_C^* E\otimes\Pc)) \\
&\iso \Rb \pi_{\widehat C\times \mathbb P^1,*}(\pi_{C\times \widehat
C}^*(\pi_C^*E\otimes\Pc))\otimes
\hat q^*(F) \\
& \iso \hat p^*(\Rb \pi_{\widehat C,*}(\pi_C^*(E)\otimes \Pc))\otimes \hat q^*(F)=
\hat p^*(\mathcal S(E))\otimes \hat q^*(F)
\end{aligned}
\label{FMbasechange}
\end{equation}
where the second isomorphism is due to $q\circ \pi_{C\times {\mathbb P^1}}=
\hat q\circ \pi_{\widehat C\times {\mathbb P^1}}$ and $p\circ \pi_{C\times
{\mathbb P^1}}= \pi_C\circ \pi_{\widehat C\times \widehat C}$, the third is
the projection formula and the forth is base change in the derived category. We
also see that given morphisms $\Phi\colon E_2\to E_1$ of vector bundles on $C$
and $\gamma\colon F_2\to F_1$ of vector bundles on $\mathbb P^1$, then the
morphism $\mathcal S_{\mathbb P^1}(p^*\Phi\otimes q^*\gamma)$ is identified
with $\hat p^*(\mathcal S(\Phi))\otimes \hat q^*\gamma$, that is, the following
diagram is commutative
\begin{equation}
\xymatrix{ {\mathcal S}_{\mathbb P^1}(p^* E_2\otimes q^* F_2)
\ar[rr]^{\mathcal S_{\mathbb P^1} (p^*\Phi\otimes q^*\gamma)} \ar[d]_{\wr}&&
{\mathcal S}_{\mathbb P^1}(p^* E_1\otimes q^* F_1)
\ar[d]^{\wr} \\
\hat p^*(\mathcal S(E_2))\otimes \hat q^*(F_2)\ar[rr]^{\hat p^*({\mathcal
S}(\Phi))\otimes \hat q^*\gamma} & & \hat p^*(\mathcal S(E_1))\otimes \hat
q^*(F_1) } \label{basechange}
\end{equation}
where the vertical isomorphisms are the base change isomorphisms \eqref{FMbasechange} we have just considered.

We shall give two natural definitions of the Fourier-Mukai transform of a
triple and show that they are equivalent under the isomorphism given in
$(ii)$ of  Proposition \ref{invv}. First we must ensure that the transform of a
triple is again a triple.

\begin{defin} The triple $T=(E_1, E_2, \Phi)$ is $\IT_i$ if both
bundles $E_1$, $E_2$ are $\IT_i$ with the same index $i$.
\end{defin}

\begin{defin1}\label{defintrans1} Let $T=(E_1, E_2, \Phi)$ be an
$\IT_i$ triple. The Fourier-Mukai transform of $T$ is defined as the triple
$\widehat T = (\mathcal S^i (E_1), \mathcal S^i   (E_2), \mathcal S^i
(\Phi))$. We shall write $\widehat T
= (\widehat E_1, \widehat E_2, \widehat \Phi)$ for the transformed triple.
\end{defin1}

Since a triple $T$ corresponds exactly to an $SU(2)$-equivariant bundle $E$ on
$C\times \mathbb P^1$, this suggests another definition of the Fourier-Mukai
transform of an $\IT_i$  triple as the triple associated to the transform of
the bundle $E$ with respect to the relative transform $\mathcal S_{\mathbb
P^1} $. This observation leads us in a natural way to consider a relative
version of the Nahm transform, an argument that will be pursued in the next
section. Note that in order that the relative transform  of the bundle $E$
consists of a single sheaf, we should ensure that
 $E$ is $\IT_i$.  This is achieved by the following Proposition whose proof is a
straightforward consequence of the base change property of the Fourier-Mukai
transform.

\begin{prop}\label{prop:extension} If $E_1$ and $E_2$ are
$\IT_i$-bundles with respect to $\mathcal S$ (with the same index $i$), then
$E$ is $\IT_i$ with respect to $\mathcal S_{\mathbb P^1}$ and its transform
$\widehat E$ sits in an exact sequence of the type $$ 0\to \hat p^* \widehat
E_1 \to \widehat E \to \hat p^*\widehat E_2 \otimes \hat q^* \mathcal
O_{\mathbb P^1} (2)\to 0.$$ \noindent Therefore $\widehat E$ is an
$SU(2)$-equivariant bundle on $\widehat C\times \mathbb P^1$.
\qed\end{prop}

Now we can define.

\begin{defin2}\label{defintrans2} We define the Fourier-Mukai
transform of a $\IT_i$ triple  $T=(E_1,E_2,\Phi)$ as the triple associated to
the transform $\widehat E= \mathcal S_{\mathbb P^1}^i (E)$ of the associated
$SU(2)$-equivariant  and $\IT_i$ bundle $E$ on $C\times \mathbb P^1$.
\end{defin2}

It remains to check that  definitions \ref{defintrans1} and
\ref{defintrans2} are compatible.

\begin{prop}\label{compa} Let $T$ be an $\IT$ triple and let
$E$ be the corresponding invariant bundle on $C\times \mathbb P^1$, then the
(absolute) Fourier-Mukai transform $\widehat T$ in Definition \ref{defintrans1}
corresponds to the triple given by the transform $\widehat E$ of Definition
\ref{defintrans2} under the isomorphism given in Proposition \ref{invv}. In
other words we have the following commutative diagram
\begin{equation*}
\xymatrix{{} \Ext^1_{C\times \mathbb P^1} (p^* E_2 \otimes q^* \mathcal
O_{\mathbb P^1}(2), p^* E_1)) \ar[r]^{\mathcal S_{\mathbb P^1}}
\ar[d]_{\sigma_C\ \wr} &
 \Ext^1_{C\times \mathbb P^1} (\hat p^* \widehat E_2
\otimes \hat q^*\mathcal O_{\mathbb P^1} (2), \hat p^*\widehat E_1))
\ar[d]^{\wr\ \sigma_{\widehat C}}
\\
\Ext^0_C (E_2, E_1) \ar[r]^{\mathcal S} & \Ext^0_C (\widehat E_2,
\widehat E_1) }
\end{equation*}
where the vertical rows are the isomorphisms introduced in Proposition
\ref{invv} and the horizontal isomorphisms are induced from the relative and
absolute Fourier-Mukai transforms.
\end{prop}
\begin{proof} Given a morphism $\Phi\colon E_2\to E_1$, we know by Remark
 \ref{explicitPhi}, that $\sigma_C^{-1}(\Phi)$ is the element of
 $\Ext^1_{C\times \mathbb P^1} (p^*E_2 \otimes q^*\mathcal O_{\mathbb P^1} (2), p^*E_1)$
 corresponding by \eqref{ext1} to the morphism
$$ p^*(\Phi) \otimes q^*(\tr^{-1})\colon p^*E_2 \otimes q^*\mathcal O_{\mathbb
P^1} (2) \to p^*E_1[1]\,. $$ Now, by \eqref{basechange}, $\mathcal S_{\mathbb
P^1}(\sigma_C^{-1}(\Phi))$ is the element of $\Ext^1_{\widehat C\times \mathbb
P^1} (\hat p^*\mathcal S(E_2) \otimes \hat q^*\mathcal O_{\mathbb P^1} (2),
\hat p^*\mathcal S(E_1))$ corresponding by \eqref{ext1} to the morphism $$
\hat p^*(\mathcal S(\Phi)) \otimes \hat q^*(\tr^{-1})\colon \hat p^*\mathcal
S(E_2) \otimes \hat q^*\mathcal O_{\mathbb P^1} (2) \to \hat p^*\mathcal
S(E_1)[1]$$ which, again by Remark \ref{explicitPhi}, corresponds to
$\sigma_{\widehat C}^{-1}(\mathcal S(\Phi))$.
\end{proof}

\begin{remark}\label{comparem}  In order to ensure that the Fourier-Mukai
transform gives rise to morphisms between moduli spaces of triples one should
check that the transform preserves families of $(\IT)$ triples. This can be
checked directly as in the usual case of families of sheaves, alternatively
one can use Remark \ref{family} and note that the Fourier-Mukai transform is
well-behaved with respect to families and therefore induces morphisms
between the moduli spaces of $SU(2)$-equivariant sheaves.
\end{remark}

\subsection{Preservation of stability for small $\alpha$}\label{preservation-s}

Let $\mathcal N_{\alpha_m^+}^s (n_1, n_2, d_1, d_2)$ be the moduli
space of  $\alpha_m^+ $-stable triples with
$\alpha_m^+=\alpha_m+\epsilon $ such that $\epsilon
>0$ and $(\alpha_m, \alpha_m^+]$ does not contain any critical value.
(We assume that $d_1/n_1 \geq d_2/n_2$, since this is a necessary
condition for the moduli space not to be empty.) One has the
following (Proposition 3.23 in \cite{BG}).

\begin{prop}\label{stabsmall} If a triple $T= (E_1, E_2,\Phi)$ is
$\alpha_m^+$-stable, $E_1$ and $E_2$ are semistable. Conversely,
if $E_1$ and $E_2$ are stable and $\Phi\neq 0$ then $T=(E_1,E_2,
\Phi)$ is $\alpha_m^+$-stable.
\qed
\end{prop}

\begin{prop}\label{modsmall} If $(n_1,d_1)=1$, $(n_2,d_2)=1$ and $d_1/n_1 > d_2/n_2$, the
moduli space of stable triples $\mathcal N_{\alpha^+_m}^s$ is isomorphic to a
$\mathbb P^N$-fibration over $\mathcal{M}_C(n_1,d_1)\times
\mathcal{M}_C(n_2,d_2)$, where $N=n_2d_1-n_1d_2-1$.
\end{prop}
\begin{proof}
By Proposition \ref{stabsmall}, $\mathcal N_{\alpha^+_m}^s$ is the
projectivization of a Picard sheaf on $\mathcal{M}_C(n_1,d_1)\times
\mathcal{M}_C(n_2,d_2)$ (Corollary 6.2 in \cite{BGG1}), which in this case is
a vector bundle with fibre $ H^0 (C, E_2^*\otimes E_1)$ over $(E_1,E_2)$,
since $ H^1 (C, E_2^\vee\otimes E_1)\simeq  H^0 (C, E_1^\vee\otimes E_2)^*=0$.
\end{proof}

Given an $\IT_i$ triple $T=(E_1, E_2, \Phi)$ with transform $\widehat T =
(\widehat E_1, \widehat E_2, \widehat \Phi)$ we denote by $\hat\alpha_m$ the
minimum value of the stability parameter $\hat\alpha$ with the type $(\hat
n_1,\hat n_2,\hat d_1,\hat d_2)$ defined by $\widehat T$. As above,
$\hat\alpha_m^+$ is any real number such that the interval $(\hat\alpha_m,
\hat\alpha_m^+]$ does not contain critical values.

\begin{thm}\label{pres2} Let $T=(E_1, E_2, \Phi)$ be a $\alpha_m^+$-stable
triple such that $(n_1,d_1)=1$, $(n_2,d_2)=1$ and $d_1d_2>0$ (this forces
$\Phi\neq 0$). Then the Fourier-Mukai transform $\widehat T = (\widehat E_1,
\widehat E_2, \widehat \Phi)$ is $\hat \alpha_m^+$-stable. The result also
holds in the converse direction with the obvious modifications on the
hypotheses.
\end{thm}

\begin{proof} By Proposition \ref{stabsmall} we have that $E_1$
and $E_2$ are both semistable. Moreover, $E_1$ and $E_2$ are stable due to the
conditions on the rank and degree. Thus in the triple $\widehat T =(\widehat
E_1, \widehat E_2, \widehat \Phi)$ both bundles are stable. By  Proposition
\ref{stabsmall} again we conclude that the triple $\widehat T$ is also
$\alpha^+_m$-stable. The proof of the converse is identical.
\end{proof}

\begin{corol}\label{prescol2} Keeping the conditions stated in the previous
Theorem and assuming additionally that $d_1/n_1>d_2/n_2$, then the
Fourier-Mukai transform induces an isomorphism $$ \mathcal S\colon \mathcal
N_{\alpha_m^+}^s \iso \mathcal N_{\hat \alpha_m^+}^s .$$ In
other words, the Fourier-Mukai transform induces an isomorphism between the
$\mathbb P^N$-fibrations described in Proposition \ref{modsmall}. \qed
\end{corol}

\subsection{Preservation of stability for large $\alpha$}\label{preservation-l}

Recall that $\alpha_L $ is the largest critical value in the
interval $(\alpha_m, \alpha_M)$. If $\alpha_L < \alpha  <\alpha_M
$ the stability condition does not vary in this range, and we can
then denote by  $\mathcal N_{\alpha_M^-}^s (n_1, n_2, d_1, d_2)$
the moduli space of  $\alpha $-stable triples for any value
$\alpha \in (\alpha_L, \alpha_M)$.

The relationship between the stability of the triple and that of
the involved  bundles is given by the following Proposition
(\cite[Propositions 7.5 and 7.6]{BGG1}).

\begin{prop}\label{form1} Let $T= (E_1, E_2, \Phi)$ be an $\alpha $-semistable
triple for some $\alpha_L < \alpha  <\alpha_M $, and let us
suppose that $n_1>n_2$. Then $T$ defines an extension of the form
\begin{equation}\label{ext}
0\to E_2 \xrightarrow{\Phi} E_1 \to F \to 0
\end{equation}
with $F$ locally free, and $E_2$ and $F$ are semistable.
Conversely, let  $T= (E_1, E_2, \Phi)$ be a triple defined by a non trivial
extension  of the form (\ref{ext}),
with $F$ locally free. If $E_2$ and $F$ are stable then $T$ is
$\alpha $-stable for $\alpha_L < \alpha <\alpha_M $.
\qed\end{prop}

{F}rom this we have the following result (Theorem 7.7 in \cite{BGG1}).

\begin{thm}\label{Bir1} Let  $n_1>n_2$, $d_1/n_1>d_2/n_2$,
$(n_1-n_2,d_1-d_2)=1$ and $(n_2,d_2)=1$. Then the moduli space $\mathcal
N_{\alpha_M^-}^s (n_1, n_2, d_1, d_2)$ is smooth of dimension $n_2d_1-n_1d_2
+1$ and it is isomorphic to a $\mathbb P^N$-fibration over $\mathcal
M_C(n_2, d_2) \times \mathcal M_C(n_1 - n_2, d_1 - d_2)$, whose  fibre
over the point $(E_2,F)$ is given by $\mathbb P H^1 (C, E_2\otimes F^*)$,
and $N=n_2d_1-n_1d_2-1$.
\end{thm}

\begin{remark} The case $n_1<n_2$  reduces to the situation in Theorem
\ref{Bir1} by considering the dual triple.
\end{remark}

We prove now that the Fourier-Mukai transform preserves stability for ``large''
values of the parameter $\alpha$.

\begin{thm}\label{pres1} Let $T=(E_1, E_2, \Phi)$ be an $\alpha$-stable
triple such that $(n_1-n_2, d_1-d_2)=1$, $(n_2, d_2)=1$, $n_1 \neq n_2$ and
$\alpha_L < \alpha <\alpha_M $. Suppose also that $d_1>0$, $d_2>0$ and $d_1 -
d_2>0$ (resp. $d_1 < 0$, $d_2<0$ and $d_1-d_2 < 0$); then $T$ is $\IT_0$ (resp.
$\IT_1$) and the transformed triple $\widehat T = (\widehat E_1, \widehat E_2,
\widehat \Phi)$ is $\hat \alpha $-stable for $\hat \alpha \in (\hat \alpha_L,
\hat \alpha_M)$ where $\hat \alpha_L$ and $\hat \alpha_M$ are the values
corresponding to the transformed triple $\widehat T$.
\end{thm}

\begin{proof} We prove the $\IT_0$ case, the proof of the other case is
entirely similar. Without loss of generality we may assume $n_1>n_2$. By
Proposition \ref{form1} the map $\Phi\colon E_2 \to E_1$ is injective and the
quotient sheaf $F$ in $0\to E_2 \to E_1 \to F \to 0$
is locally free. Moreover, $E_2$ and $F$ are stable, and hence
$\IT_0$, from which  it follows that $E_1$ is $\IT_0$. Transforming the above sequence we
get
\begin{equation*}
0\to \widehat E_2 \to \widehat E_1 \to \widehat F \to 0\,.
\end{equation*}
Since the Fourier-Mukai transform preserves stability (Proposition
\ref{presstab}) it follows that $\widehat E_2$ and $\widehat F$ are stable.
By Proposition \ref{form1} $\widehat T = (\widehat E_1,
\widehat E_2, \widehat \Phi)$ is $\hat \alpha$-stable for $\hat \alpha \in
(\hat \alpha_L, \hat \alpha_M) $. The proof of the converse is identical.
\end{proof}

Under the same conditions of Theorem \ref{pres1} we have the following.

\begin{corol}\label{prescol1}  The Fourier-Mukai transform induces an
isomorphism between the moduli spaces of $\IT_i$ stable triples:
$$ \mathcal N_{\alpha_M^-}^s (n_1,n_2, d_1, d_2) \simeq \mathcal
N_{\hat \alpha_M^-}^s ((-1)^i d_1, (-1)^i d_2, (-1)^{i+1} n_1,
(-1)^{i+1} n_2)\,.$$ As a consequence, the Fourier-Mukai transform
yields an isomorphism between the $\mathbb{P}^N$-fibrations
described in Theorem \ref{Bir1}.
\end{corol}

\subsection{Applications to moduli spaces on $C\times \mathbb P^1$}

One notable application of the theory of triples is the
construction of slope-stable bundles on $C\times \mathbb P^1$ with
respect to the polarization $\omega_\alpha$, with $\alpha>0$ (see
Theorem 9.2 in \cite{BGG1}). It seems quite natural to use the
relative transform $\mathcal S_{\mathbb P^1}$ to further study the
properties of those bundles and to produce new examples of stable
bundles. We give in this section a result on the preservation of
stability for a class of bundles on $C\times \mathbb P^1$ which
can not be handled using the standard techniques based on choosing
``suitable polarizations'' as done for example in \cite{Br1} or
\cite{HP}, because the polarizations $\omega_\alpha$ are not
suitable; the reason for this being that there exist
$SU(2)$-equivariant bundles which are $\omega_\alpha$-stable and
whose restriction to a fibre, is never stable (here we are
assuming $\Phi \neq 0$). To see this, take $E$ such that
$\Phi\colon E_2\to E_1$ is not an isomorphism and note that the
restriction of such a bundle to a fibre $C_t$ is given by an
extension $$ 0\to E_1 \to E_t \to E_2 \to 0\,. $$ Since the
associated triple is stable, Lemma 4.5 in \cite{BG} implies that
$\mathrm{Ext}^1(E_2, E_1) =0$ whenever $\Phi$ is not an isomorphism, therefore the
previous extension is always split and the restriction $E_t$ is
not stable.

The following Proposition follows now immediately.

\begin{prop}
Let $T=(E_1,E_2,\Phi)$ be an $\alpha$-stable triple and let $E$ be its
associated vector bundle on $C\times\mathbb{P}^1$. Then $T$ is $\IT_i$ if and
only if $E$ is $\IT_i$ with respect to $\mathcal{S}_{\mathbb{P}^1}$.
\qed
\end{prop}

We can use this Proposition to prove the following result.

\begin{thm}\label{suinvone} Let $T=(E_1,E_2,\Phi)$ be an $\alpha$-stable triple
with $E_1\simeq E_2$ and $\Phi\neq 0$. Assume that either $\rk(E_1)=\rk(E_2)>1$ or
$\deg E_1=\deg E_2\neq 0$. Then the associated $SU(2)$-equivariant
bundle $E$ on $C\times\mathbb{P}^1$ is IT and the Fourier-Mukai transform
$\widehat E$ is polystable. Moreover, the triple $\widehat T=(\widehat
E_1,\widehat E_2,\widehat \Phi)$ is $\hat\alpha$-stable for any $\hat\alpha>0$.
\end{thm}
\begin{proof} Thanks to Theorem \ref{dimred} we have $ E\simeq
(p^* E_1 \otimes q^* \mathcal O_{\mathbb P^1} (1) ) \oplus ( p^* E_2 \otimes
q^* \mathcal O_{\mathbb P^1} (1) ) $ with $E_1\simeq E_2$ stable. The base
change property for the Fourier-Mukai transform implies $$
\mathcal{S}_{\mathbb{P}^1}(E)\simeq (p^* \mathcal{S}(E_1) \otimes q^* \mathcal
O_{\mathbb P^1} (1) ) \oplus ( p^* \mathcal{S}(E_2) \otimes q^* \mathcal
O_{\mathbb P^1} (1) ) $$ Therefore, $E$ is $\IT$ if and only if $E_1\simeq
E_2$ is $\IT$ with respect to $\mathcal{S}$, and this follows from Proposition
\ref{presstab}  since the stability of $E_1\simeq E_2$ implies that its degree
is not zero unless the rank is 1. The polystability of $\widehat E$ is a
consequence of the above expression for $\widehat E$ and the fact that
$\mathcal{S}$ preserves stability, see Proposition \ref{presstab}.

On the other hand, let us recall that a triple $(E_1,E_2,\Phi)$ with $E_1\simeq
E_2$ is $\alpha$-stable, for any $\alpha>0$, if and only if $\Phi$ is an
isomorphism and $E_1\simeq E_2$ is stable \cite[Lemma 4.6]{BG}. These
conditions are preserved by $\mathcal{S}$, therefore $\widehat T$ is
$\hat\alpha$-stable for any $\hat\alpha>0$.
\end{proof}

Collecting previous results, particularly  Theorem \ref{pres1}, Theorem
\ref{pres2} and Theorem \ref{dimred}, we have.

\begin{thm}\label{suinv} Let $T$ be an $\alpha$-stable triple of type
$(n_1, n_2, d_1, d_2)$ with $\alpha$ and $(n_i,d_i)$ satisfying one of the
conditions
\begin{itemize}
\item[(i)] $(n_1-n_2, d_1-d_2)=1$, $(n_2,d_2)=1$, $n_1 \neq n_2$ and $\alpha_L <
\alpha <\alpha_M $. Suppose also that $d_i>0$ for $i=1,2$, $d_1 - d_2>0$ (resp.
$d_i < 0$ $i=1,2$, $d_1-d_2 < 0$) and $\alpha_L < \alpha < \alpha_M$,
\item[(ii)] $(n_1,d_1)=1$, $(n_2,d_2)=1$, $d_1d_2>0$ and
$\alpha_m < \alpha < \alpha_m^+$,
\end{itemize}
(i.e.  one of the conditions in Theorems \ref{pres1} or \ref{pres2}). Then,
the corresponding $SU(2)$-equivariant bundle $E$ on $C\times \mathbb P^1$ is
$\IT$. Moreover if $E_1$ and $E_2$ are not isomorphic, then the Fourier-Mukai
transform $\widehat E$ is stable with respect to the polarization
$\omega_{\hat \alpha}$, where $\hat \alpha $ is the corresponding parameter
for the transformed triple according to Theorem \ref{pres1} in case $(i)$ and
to Theorem \ref{pres2} in case $(ii)$. \qed
\end{thm}

The relative Fourier-Mukai transform induces an isomorphism
between the corresponding moduli spaces of $SU(2)$-equivariant
bundles as follows from the previous Theorem and Remarks
\ref{family} and \ref{comparem}. Therefore we have.

\begin{corol} Let $\mathcal N_\alpha^s$ be a
moduli space of $\alpha$-stable triple satisfying one of the
conditions (i) or (ii). Let $\mathcal M_\alpha^{SU(2)}$ be the
corresponding moduli space of $SU(2)$-equivariant bundles on
$C\times \mathbb P^1$. Then the relative Fourier-Mukai transform
gives an isomorphism
 $$ \mathcal S_{\mathbb
P^1} \colon \mathcal M_\alpha^{SU(2)} \iso {\mathcal
M}_{\hat\alpha}^{ SU(2)}\,. $$ \qed
\end{corol}

\section{Nahm transforms for triples}\label{transnahm}

\subsection{Relative Nahm transform}
In this section we modify the absolute Nahm transform to produce a relative version of it.

For every elliptic curve $C$ we consider the projections
$q\colon X=C\times \mathbb{P}^1 \to\mathbb{P}^1$, $\hat q\colon \widehat
X=\widehat C\times \mathbb{P}^1\to \mathbb{P}^1$ where $\widehat
C$ is the dual elliptic curve. We endow
 the pull-back $\mathcal{P}_{\mathbb{P}^1}$ of the Poincar\'{e} line
bundle to $X\times_{\mathbb{P}^1}\widehat X$, with the pull-back connection
$\nabla_{\mathcal{P}_{\mathbb{P}^1}}$. For every point $\hat
x=(\xi,t)\in\widehat C\times\mathbb{P}^1$ we endow the Hermitian
line bundle $\mathcal{P}_{{\mathbb{P}^1},\,\hat
x}\equiv{\mathcal{P}_{\mathbb{P}^1}}_{\vert{X_{\hat q(\hat
x)}}}\to X_{\hat q(\hat x)}$,  obtained by restricting
$\mathcal{P}_{\mathbb{P}^1}$ to the fiber $X_{\hat q(\hat
x)}\subset X\times_{\mathbb{P}^1}\widehat X$ of $q$ over $\hat
q(\hat x)\in\mathbb{P}^1$, with the flat unitary connection
$\overline\nabla_{\hat{x}}$ given by the restriction of
$\nabla_{\mathcal{P}_{\mathbb{P}^1}}$. In this way $\widehat X$
parametrizes the gauge equivalence classes of Hermitian flat line
bundles along the fibers of $q\colon X\to\mathbb{P}^1$.

Let us consider a Hermitian vector bundle $E\to X$ with a unitary
connection $\nabla$. We denote by $E_t$ the restriction of $E$ to
the fibre $X_t=q^{-1}(t)$, $\nabla_{t}$ is the restriction of
$\nabla$ to $E_t$. On the vector bundle $E_{\hat q(\hat
x)}\otimes\mathcal{P}_{{\mathbb{P}^1},\,\hat x}$ we have the
connection $\nabla_{\hat x}=\nabla_{\hat q(\hat x)}\otimes
1+1\otimes \overline\nabla_{\hat x}$. Therefore we have the family
of coupled Dirac operators $$D_{\hat x}=\sqrt
2\bar\partial^*_{E_{\hat q(\hat x)}\otimes
\mathcal{P}_{{\mathbb{P}^1},\,\hat x}}\colon \Omega^0(X_{\hat q(\hat
x)}, E_{\hat q(\hat x)}\otimes \mathcal{P}_{{\mathbb{P}^1},\,\hat
x})\to \Omega^{0,1}(X_{\hat q(\hat x)}, E_{\hat q(\hat x)}\otimes
\mathcal{P}_{{\mathbb{P}^1},\,\hat x}).$$

As in the absolute case we define the index $\mathrm{Ind}(D)$ of this
family of Dirac operators $D$ parametrized by $\widehat C\times \mathbb{P}^1$.
The relative Nahm transform maps a Hermitian vector bundle with a unitary
connection over $C\times\mathbb{P}^1$ into a Hermitian vector bundle with a
unitary connection over $\widehat C\times \mathbb{P}^1$.

\begin{defin}
Let $(E,\nabla)$ be a pair formed by a Hermitian vector bundle $E$
over $C\times\mathbb{P}^1$ and a unitary connection $\nabla$ on
$E$. We say that $(E,\nabla)$ is an $\IT_{\mathbb{P}^1}$ (index
Theorem) pair relative to ${\mathbb{P}^1}$ if either $\Coker D=0$
or $\Ker D=0$. In the first case we say that $(E,\nabla)$ is an
$\IT_{{\mathbb{P}^1},\,0}$-pair, whereas in the second we call it
an $\IT_{{\mathbb{P}^1},\,1}$ pair. The transformed bundle of an
$\IT_{{\mathbb{P}^1},\,i}$-pair is, according to the parity of i,
the vector bundle $\widehat E=\pm\mathrm{Ind}(D)\to \widehat
C\times\mathbb{P}^1$.
\end{defin}

Proceeding  in the same way as in the absolute case we can endow
the transformed vector bundle of an $\IT_{\mathbb{P}^1}$-pair with
a Hermitian metric and a unitary connection in a natural way. In
doing this, since all the fibrations involved are trivial, the
main difference one encounters is that the parameter space of the
family is enlarged from $\widehat C$ to $\widehat
C\times\mathbb{P}^1$, but since $X_{\hat q(\hat x)}\simeq C$ the
Dirac operators are still defined on vector bundles over the
elliptic curve $C$. Therefore, the theory parallels the one
developed in the absolute setting.

\begin{defin}
Let $(E,\nabla)$ be an $\IT_{\mathbb{P}^1}$-pair. We call $(\widehat
E,\widehat\nabla)$ the relative Nahm transform of $(E,\nabla)$ and denote it
by $\mathcal{N}_{\mathbb{P}^1}(E,\nabla)$.
\end{defin}

Let $E\to C\times \mathbb{P}^1$ be a holomorphic vector bundle
endowed with a unitary connection $\nabla$ compatible with the
holomorphic structure. Since the spin$^c$ Dirac operator $D_{\hat
x}$ gets identified with the Dolbeault-Dirac operator of $E_{\hat
q(\hat x)}\otimes \mathcal{P}_{{\mathbb{P}^1},\,\hat x}$,  by
Hodge theory and the Dolbeault isomorphism we have
\begin{align}\label{dolb1} \Ker D_{\hat
x}&\simeq H^{0}(X_{\hat q(\hat x)},E_{\hat q(\hat x)}
\otimes \mathcal{P}_{{\mathbb{P}^1},\,\hat x})\\
\label{dolb2} \Coker D_{\hat x} &\simeq H^{1}(X_{\hat q(\hat
x)},E_{\hat q(\hat x)}\otimes \mathcal{P}_{{\mathbb{P}^1},\,\hat
x}).\end{align}

Let us suppose that $E$ is $\IT_i$ with respect to the relative Fourier-Mukai
transform described in Section \ref{fm-triple}. The isomorphisms (\ref{dolb1})
and (\ref{dolb2}) mean that $(E,\nabla)$ is an
$\IT_{{\mathbb{P}^1},\,i}$-pair. As we saw there,  by \cite[Theorem 2]{BBH1}
or \cite[Theorem 3.2.8]{DK} we have a natural $C^\infty$ vector bundle
isomorphism induced by Hodge theory $$\phi_{\mathbb{P}^1}\colon \mathcal
S^i_{\mathbb{P}^1}(E)\iso\widehat E.$$ Moreover, since the
Poincar\'{e} bundle $\mathcal{P}_{\mathbb{P}^1}\to C\times\widehat C$ is a
holomorphic line bundle and the connection
$\nabla_{\mathcal{P}_{\mathbb{P}^1}}$ is compatible with the holomorphic
structure, the same arguments that in the absolute case led us to prove
Theorem \ref{compatibility} give us now the following.

\begin{thm}\label{compatibility2} Let $F_1$, $F_2$ be Hermitian vector
bundles over $C\times\mathbb{P}^1$ endowed with unitary
connections $\nabla_1$, $\nabla_2$ such that $(F_1,\nabla_1)$,
$(F_2,\nabla_2)$ are $\IT_{\mathbb{P}^1,i}$-pairs with respect to
the Nahm transform. Then we have
\begin{enumerate}\item The connections $\widehat\nabla_1$,
$\widehat\nabla_2$ are compatible with the holomorphic structures
of $\mathcal S^i_{\mathbb{P}^1} (F_1)$, $\mathcal
S^i_{\mathbb{P}^1} (F_2)$, respectively. Thus, the curvature of
the connections $\widehat\nabla_1$, $\widehat\nabla_2$ is of type
$(1,1)$.
\item For every holomorphic morphism $\Psi\colon F_1\to F_2$
we have an induced holomorphic morphism
$\mathcal{N}(\Phi)\colon\widehat F_1\to \widehat F_2$ and a
commutative diagram
\begin{equation}
\xymatrix{{\widehat
F_1}\ar[r]^{\phi_{F_1}}\ar[d]_{\mathcal{N}_{\mathbb{P}^1}(\Psi)} &
{\mathcal{S}(F_1)}\ar[d]^{\mathcal{S}_{\mathbb{P}^1}(\Psi)}\\
{\widehat F_2}\ar[r]^{\phi_{F_2}} & {\mathcal{S}(F_2)} }
\label{comm2}
\end{equation}
\end{enumerate}
\end{thm}

\subsection{Relative Nahm transform for $SU(2)$-invariant Einstein-Her\-mitian
connections}

Let us suppose that $E_1$, $E_2$ are complex Hermitian vector
bundles over $C$ and let us choose an $SU(2)$-invariant metric on
$H^{\otimes2}$. We put on the bundle $E=p^*E_1\oplus(p^*E_2\otimes
q^*H^{\otimes2})$  the Hermitian metric which is determined in a
natural way by the Hermitian metrics of $E_1$, $E_2$ and
$H^{\otimes2}$.

By Proposition 3.5 in \cite{G1} there is a one to one correspondence between
the $SU(2)$-invariant unitary connections on $E$ and the triples
$\mathcal{T}=((E_1,\nabla_1),(E_2,\nabla_2),\Phi)$ formed by unitary
connections $\nabla_1$, $\nabla_2$ on $E_1$, $E_2$, respectively, and a
$C^\infty$ vector bundle morphism $\Phi\colon E_2\to E_1$. Moreover, this
correspondence also holds at the level of $SU(2)$-invariant holomorphic
structures on $E$. Before discussing it we introduce the following.

\begin{defin} We call a triple $\mathcal{T}=((E_1,\nabla_1),(E_2,\nabla_2),\Phi)$
integrable if $\Phi\colon E_2\to E_1$ is holomorphic with respect to the
holomorphic structures determined by the connections $\nabla_1$ and $\nabla_2$.
\end{defin}

Proposition 3.9 in \cite{G1} gives us a one to one correspondence
between $SU(2)$-invariant holomorphic structures on $E$,
considered as integrable $SU(2)$-invariant connections (i.e.
connections with curvature of type $(1,1)$), and integrable
triples $\mathcal{T}=((E_1,\nabla_1),(E_2,\nabla_2),\Phi)$. This
is precisely the content of Proposition \ref{invv} which gives us
a bijective correspondence between the $SU(2)$-invariant
holomorphic structures on $E$ and holomorphic triples
$T=(\mathcal{E}_1=(E_1,\bar\partial^{\nabla_1}),
\mathcal{E}_2=(E_2,\bar\partial^{\nabla_2}),\Phi)$.

Let us denote by $\nabla^{\mathcal{T}}$ the $SU(2)$-invariant integrable
connection on $E$ determined by an integrable triple $\mathcal{T}$. If we
express its curvature with respect to the splitting
$E=p^*E_1\oplus(p^*E_2\otimes q^*H^2)$ we have

\begin{equation}\label{curvature}
F^{\nabla^{\mathcal{T}}}=\begin{pmatrix} p^*F^{\nabla_1}-\beta\wedge\beta^* &
\partial\beta\\-\bar\partial\beta^* & p^*F^{\nabla_2}\otimes
1+1\otimes q^*F^{\nabla^\prime}-\beta^*\wedge\beta,
\end{pmatrix}
\end{equation} where $F^{\nabla_i}$ is the curvature of the
connection $\nabla_i$, $F^{\nabla^\prime}$ is the curvature of the unique
$SU(2)$-invariant unitary connection on $H^{\otimes2}$, $\beta=p^*\Phi\otimes
q^*\eta$, with $\eta$ an $SU(2)$-invariant section of $H^{\otimes -2}$ and
$\bar\partial$ is the Cauchy-Riemann operator determined by the connections
$\nabla_1$, $\nabla_2$ and $\nabla^\prime$, for further details see
\cite{G1,BG}.

We want to study the relative Nahm transform of the $SU(2)$-equivariant
bundles $(E,\nabla^{\mathcal{T}})$ associated to integrable triples.

The following is straightforward.

\begin{prop} Let $\mathcal{T}=((E_1,\nabla_1),(E_2,\nabla_2),\Phi)$ be an
integrable triple on $C$ and let $(E,\nabla^{\mathcal{T}})$ be
its associated bundle with connection over $C\times
\mathbb{P}^1$. If both $(E_1,\nabla_1)$, $(E_2,\nabla_2)$ are
$\IT_i$-pairs then $(E,\nabla^{\mathcal{T}})$ is an
$\IT_{\mathbb{P}^1,\,i}$ pair.
\qed\end{prop}

Given an integrable triple $\mathcal{T}=((E_1,\nabla_1),(E_2,\nabla_2),\Phi)$
such that $(E_1,\nabla_1)$ and $(E_2,\nabla_2)$ are $\IT_i$-pairs we can form
the triple $\widehat{\mathcal{T}}=((\widehat E_1,\widehat \nabla_1),(\widehat
E_2,\widehat\nabla_2),\widehat\Phi)$ obtained by means of the absolute Nahm
transform.  Here we have denoted by $\widehat\Phi$ the Nahm transform
$\mathcal{N}(\Phi)$. By the sake of brevity the same notation is used
hereafter. On the other hand, if $(E,\nabla^{\mathcal{T}})$ is the vector
bundle with connection over $C\times \mathbb{P}^1$ associated to the triple
$\mathcal{T}$, we can apply to it the relative Nahm transform to obtain
$\mathcal{N}_{\mathbb{P}^1}(E,\nabla^{\mathcal{T}})$. Taking into account the
compatibility between the Fourier-Mukai and Nahm transforms, Theorems
\ref{compatibility} and \ref{compatibility2} and Proposition \ref{compa} we
have.

\begin{prop}
$\mathcal{N}_{\mathbb{P}^1}(E,\nabla^{\mathcal{T}})$ is the vector bundle on
$\widehat C\times \mathbb{P}^1$ associated to the triple
$\widehat{\mathcal{T}}=((\widehat E_1,\widehat \nabla_1),(\widehat
E_2,\widehat\nabla_2),\widehat\Phi)$.
\qed\end{prop}

\begin{defin}
Let $\mathcal{T}=((E_1,\nabla_1),(E_2,\nabla_2),\Phi)$ be an integrable triple
on $C$. We say that it satisfies the $\tau$-coupled vortex equations if
\begin{align*} i\Lambda F^{\nabla_1}+\Phi\Phi^* &=2\pi\tau
\Id_{E_1}\\ i\Lambda F^{\nabla_2}-\Phi^*\Phi &=2\pi\tau^\prime \Id_{E_2},
\end{align*}
\end{defin}
Note that in order to have solutions $\tau$, $\tau^\prime$  must fulfill the
following equation
\begin{equation}\label{tau}
n_1\tau+n_2\tau^\prime=d_1+d_2,
\end{equation} with $n_i=\mathrm{rank}(E_i)$ and
$d_i=\mathrm{deg}(E_i)$.

The following Proposition was proved in \cite{G} (see also \cite{BG}).

\begin{prop}
Let $\mathcal{T}=((E_1,\nabla_1),(E_2,\nabla_2),\Phi)$ be an
integrable triple and let $\nabla^{\mathcal{T}}$ be the
corresponding connection on $E$. Let $\tau$ and $\tau^\prime$ be
related by (\ref{tau}) and let us suppose that $$
\alpha=\frac{(n_1+n_2)\tau-d_1-d_2}{n_2}>0 .$$ Then
$\mathcal{T}=((E_1,\nabla_1),(E_2,\nabla_2),\Phi)$ satisfies the
$\tau$-coupled vortex equations if and only if
$\nabla^{\mathcal{T}}$ is a Einstein-Hermitian connection on $E\to
C\times \mathbb{P}^1$ with respect to the K\"ahler form
$\omega_\alpha=\frac{\alpha}{2}p^*\omega_C+q^*\omega_{\mathbb{P}^1}
$, where $\omega_{\mathbb{P}^1}$ is the Fubini-Study K\"ahler form
normalized to volume one and $\omega_C$ is a K\"ahler form of unit
volume.
\qed\end{prop}

\begin{prop}\label{vanishing5}
Let $\mathcal{T}=((E_1,\nabla_1),(E_2,\nabla_2),\Phi)$ be an
integrable triple on $C$ which satisfies the $\tau$-coupled vortex
equations and let $X=C\times\mathbb{P}^1$. Then:

\begin{enumerate}
\item[(i)] If the Hermitian endomorphisms $2\pi\tau \Id_{E_1}-\Phi\Phi^*$ and
$2\pi\tau^\prime \Id_{E_2}+\Phi^*\Phi$ are non-negative and there exist
$x_1,\,x_2\in C$ such that $2\pi\tau \Id_{E_1}-\Phi\Phi^*(x_1)>0$,
$2\pi\tau^\prime \Id_{E_2}+\Phi^*\Phi(x_2)>0$, then $(E,\nabla^{\mathcal{T}})$
is an $\IT_{\mathbb{P}^1,\,0}$-pair and $(E_1,\nabla_1)$, $(E_2,\nabla_2)$ are
$\IT_0$ pairs.

\item[(ii)] If the Hermitian endomorphisms $2\pi\tau \Id_{E_1}-\Phi\Phi^*$ and
$2\pi\tau^\prime \Id_{E_2}+\Phi^*\Phi$ are non-positive and there exist
$x_1,\,x_2\in C$ such that $2\pi\tau \Id_{E_1}-\Phi\Phi^*(x_1)<0$,
$\pi\tau^\prime \Id_{E_2}+\Phi^*\Phi(x_2)<0$, then $(E,\nabla^{\mathcal{T}})$
is an $\IT_{\mathbb{P}^1,\,1}$-pair and $(E_1,\nabla_1)$, $(E_2,\nabla_2)$ are
$\IT_1$ pairs.
\end{enumerate}
\end{prop}

\begin{proof}
For every $\hat x=(\xi,t) $ the restriction of
$E=p^*E_1\oplus(p^*E_2\otimes q^*H^{\otimes 2})$ to $X_{\hat q(\hat x)}\simeq
C$ is isomorphic to $E_1\oplus E_2$ as $C^\infty$ bundles. Now
(\ref{curvature}) implies that the curvature of $\nabla_{\hat q(\hat x)}$ with
respect to the splitting $E_{\hat q(\hat x)}\simeq E_1\oplus E_2$ is $$
F^{\nabla_{\hat q(\hat x)}}=\begin{pmatrix} F^{\nabla_1} & 0\\0 & F^{\nabla_2}
\end{pmatrix}.
$$ The claim now follows from Theorem \ref{vanishing2}.
\end{proof}

\subsection{Covariantly constant triples}

\begin{defin} Let $\mathcal{T}=((E_1,\nabla_1),(E_2,\nabla_2),\Phi)$ be an
integrable triple on $C$. We will say that $\mathcal{T}$ is covariantly
constant if $\Phi\Phi^*$ is covariantly constant with respect to $\nabla_1$
and $\Phi^*\Phi$ is covariantly constant with respect to $\nabla_2$.
\end{defin}

\begin{rem} Denote by $\nabla$ the connection naturally induced on
$\mathrm{Hom}(E_2,E_1)$ by $\nabla_1$ and $\nabla_2$. If $\Phi$ is covariantly
constant with respect to $\nabla$ then it is easy to check that $\mathcal{T}$
is covariantly constant. Moreover, $\Phi$ is covariantly constant with respect
to $\nabla$ if and only if $\Phi\colon E_2\to E_1$ is an anti-holomorphic map.
\end{rem}

\begin{prop}\label{decomp}
Let $\mathcal{T}=((E_1,\nabla_1),(E_2,\nabla_2),\Phi)$ be a
covariantly constant integrable triple on $C$. Then we have
holomorphic orthogonal decompositions
\begin{align*}
E_1 &\simeq\Ker \Phi^* \oplus E_1^\prime\\
E_2 &\simeq\Ker \Phi^{\ } \oplus E_2^\prime
\end{align*}which are compatible with the connections, and $\Phi$ induces an
holomorphic isomorphism $\Phi\colon E_2^\prime\to E_1^\prime$.
\end{prop}

\begin{proof}  Since $\Phi\Phi^*$ and $\Phi^*\Phi$ are covariantly
constant vector bundle endomorphisms, they are holomorphic and their eigenvalues
are constant. Moreover, $\Phi\Phi^*$, $\Phi^*\Phi$ are positive Hermitian
endomorphisms whose spectrum may differ only at $0$; therefore we have
orthogonal decompositions
\begin{align*} E_1 & =\Ker \Phi^*\oplus E_1(\lambda_1)\oplus\cdots\oplus
E_1(\lambda_k)\\
E_2 & =\Ker \Phi^{\ }\oplus E_2(\lambda_1)\oplus\cdots\oplus
E_2(\lambda_k),
\end{align*}
where $E_1(\lambda_i)$, $E_2(\lambda_i)$ are the eigenbundles with eigenvalue
$\lambda_i\neq 0$ with respect to the holomorphic endomorphisms $\Phi\Phi^*$
and $\Phi^*\Phi$, respectively. Since these endomorphisms are covariantly
constant, the subbundles $E_1(\lambda_i)$, $E_2(\lambda_i)$ are preserved by
the connections $\nabla_1$, $\nabla_2$, respectively. Moreover, for every
$\lambda_i$ we have an isomorphism $$\Phi\colon
E_2(\lambda_i)\iso E_1(\lambda_i)$$ Therefore if we denote
$E_1^\prime=E_1(\lambda_1)\oplus\cdots\oplus E_1(\lambda_k)$,
$E_2^\prime=E_2(\lambda_1)\oplus\cdots\oplus E_2(\lambda_k)$, we have an
isomorphism $$\Phi\colon E_2^\prime\iso E_1^\prime$$ as required.

\end{proof}
With the same notations as above we have the following

\begin{prop}\label{cov-const}
Let $\mathcal{T}=((E_1,\nabla_1),(E_2,\nabla_2),\Phi)$ be a covariantly
constant integrable triple on $C$. Then $\mathcal{T}$ satisfies the
$\tau$-coupled equations if and only if
\begin{enumerate}
\item $\nabla_1$ induces a constant central curvature  connection on
$\Ker \Phi^*$  with factor $2\pi\tau$, unless $\Ker \Phi^*= 0$, and a
constant central curvature  connection on $E_1^\prime$ with factor
$\pi(\tau+\tau^\prime)$ unless $E_1^\prime= 0$.

\item $\nabla_2$ induces a constant central curvature  connection on
$\Ker \Phi$  with factor $2\pi\tau^\prime$ unless $\Ker \Phi=0$ and a
constant central curvature  connection on $E_2^\prime$ with factor
$\pi(\tau+\tau^\prime)$ unless $E_2^\prime= 0$.
\end{enumerate}
\end{prop}

\begin{proof}
Since $\mathcal{T}$ is covariantly constant we have the decompositions
\begin{align*} E_1 & =\Ker \Phi^*\oplus E_1(\lambda_1)\oplus\cdots\oplus
E_1(\lambda_k)\\
E_2 & =\Ker \Phi^{\ }\oplus E_2(\lambda_1)\oplus\cdots\oplus
E_2(\lambda_k),\end{align*} provided by Proposition \ref{decomp}. Moreover,
since $\mathcal{T}$ satisfies the $\tau$-coupled equations we have
\begin{align*} i\Lambda F^{\nabla_1} &=2\pi\tau
\Id_{E_1}-\Phi\Phi^*\\ i\Lambda F^{\nabla_2} &=2\pi\tau^\prime \Id_{E_2}
+\Phi^*\Phi.
\end{align*}

Therefore we have \begin{align} {i\Lambda F^{\nabla_1}}_{|\Ker \Phi^*}
&=2\pi\tau \Id_{\Ker \Phi^*} & {i\Lambda F^{\nabla_1}}_{|E_1(\lambda_i)}
&=(2\pi\tau-\lambda_i)
\Id_{E_1(\lambda_i)}\label{EH1}\\
{i\Lambda F^{\nabla_2}}_{|\Ker \Phi} & =2\pi\tau^\prime \Id_{\Ker \Phi} &
{i\Lambda F^{\nabla_2}}_{|E_2(\lambda_i)} &=(2\pi\tau^\prime+\lambda_i)
\Id_{E_2(\lambda_i)}.\label{EH2}
\end{align}

This implies that $\Ker \Phi^*$, $E_1(\lambda_i)$, $\Ker \Phi$,
$E_2(\lambda_i)$ are bundles with constant central curvature connection
with slopes
\begin{align*}
\mu(\Ker \Phi^*)& =\tau  & \mu(E_1(\lambda_i)) & =\tau-\frac{\lambda_i}{2\pi}\\
\mu(\Ker \Phi) & =\tau^\prime & \mu(E_2(\lambda_i))
& =\tau^\prime+\frac{\lambda_i}{2\pi}.
\end{align*}
But since $E_1(\lambda_i)$ is isomorphic to
$E_2(\lambda_i)$ we must have $\mu(E_1(\lambda_i))=\mu(E_2(\lambda_i))$, that
is $\lambda_i=\pi(\tau-\tau^\prime)$. Therefore in the above decompositions
there is only one eigenvalue and if we substitute $\lambda_i$ in
\ref{EH1}-\ref{EH2} we get the required values for the factors of the
constant central curvature connections.
The converse statement is just a simple checking.
\end{proof}

\begin{corol} Let $\mathcal{T}=((E_1,\nabla_1),(E_2,\nabla_2),\Phi)$ be a covariantly
constant integrable triple on $C$ which satisfies the
$\tau$-coupled equations. Then $\mathcal{T}$ is
$(\tau-\tau^\prime)$-polystable. Moreover, $\mathcal{T}$
decomposes as a sum of $(\tau-\tau^\prime)$-polystable triples.
$$\mathcal{T}=(\Ker \Phi^*,0,0)\oplus(E_1^\prime,E_2^\prime,\Phi)
\oplus(0,\Ker \Phi,0).$$
\end{corol}

\begin{proof} The Hitchin-Kobayashi correspondence
for triples, \cite[Theorem 5.1]{BG}, establishes an equivalence between
triples that satisfy the $\tau$-coupled vortex equations and
$(\tau-\tau^\prime)$-polystable triples. Therefore, the Corollary follows at
once.

However, in the present case it is possible to give a direct proof. Since
$\mathcal{T}$ satisfies the $\tau$-coupled vortex equations, Proposition
\ref{cov-const} implies that $\Ker \Phi^*$, $E_1^\prime\simeq E_2^\prime$ and
$\Ker \Phi$ are polystable bundles with slopes $$ \mu(\Ker \Phi^*) =\tau,\qquad
\mu(E_1^\prime)=\mu(E_2^\prime) =\frac{1}{2}(\tau+\tau^\prime),\qquad \mu(\Ker
\Phi) =\tau^\prime \,. $$ Therefore we have
\begin{align*}
\mu_\alpha(\Ker \Phi^*,0,0) &=\mu(\Ker \Phi^*)=\tau\\
\mu_\alpha(E_1^\prime,E_2^\prime,\Phi)=\mu(E_2^\prime) &=\mu(E_1^\prime)+
\frac{\alpha}{2}=\tau\\
\mu_\alpha(0,\Ker \Phi,0) &=\mu(\Ker \Phi)+\alpha=\tau\,,
\end{align*}
where $\alpha=(\tau-\tau^\prime)$. Since $E_2^\prime$ carries a constant
central curvature, there exists an orthogonal decomposition $$E_2^\prime  =
E_2^{(1)}\oplus\cdots\oplus E_2^{(m)}$$ compatible with the connection and
such that every factor carries an irreducible constant central curvature
connection. Since $\Phi^*\Phi=\lambda\,\mathrm{Id}_{E_2}$ it follows that we
have an orthogonal decomposition $$E_1^\prime =
\Phi(E_2^{(1)})\oplus\cdots\oplus \Phi(E_2^{(m)})$$ Thus, the triple
$(E_1^\prime,E_2^\prime,\Phi_{|E^\prime_2})$ splits into the direct sum of
subtriples $(E_2^{(i)},\Phi(E_2^{(i)}),\Phi_{\vert E_2^{(i)}})$ with
$E_2^{(i)}$ stable and $\Phi_{\vert E_2^{(i)}}$ an isomorphism. By
\cite[Proposition 3.21]{BG} this implies that
$(E_1^\prime,E_2^\prime,\Phi_{|E^\prime_2})$ is $\alpha$-polystable.
Therefore, $\mathcal{T}$ is $\alpha$-polystable (see \cite[Definition
3.15]{BG}).
\end{proof}

\begin{remark}\label{isocov}
 If $\mathcal{T}=((E_1,\nabla_1),(E_2,\nabla_2),\Phi)$ is a covariantly
constant integrable triple on $C$ which is $\alpha$-stable with $E_1\neq0$
and $E_2\neq 0$, then the previous Corollary implies that $\Phi$ has to be an isomorphism.
\end{remark}

As a consequence of Proposition \ref{vanishing5} we immediately obtain.

\begin{lemma}
Let $\mathcal{T}=((E_1,\nabla_1),(E_2,\nabla_2),\Phi)$ be a covariantly
constant integrable triple on $C$ which satisfies the $\tau$-coupled equations.
\begin{enumerate}
\item[(i)] If $\tau>0$ and $\tau^\prime>0$ then $(E,\nabla^{\mathcal{T}})$ is an
$\IT_{\mathbb{P}^1,\,0}$-pair and $(E_1,\nabla_1)$, $(E_2,\nabla_2)$ are
$\IT_0$ pairs.

\item[(ii)] If $\tau<0$ and $\tau^\prime<0$ then $(E,\nabla^{\mathcal{T}})$ is an
$\IT_{\mathbb{P}^1,\,1}$-pair and $(E_1,\nabla_1)$, $(E_2,\nabla_2)$ are
$\IT_1$ pairs.
\end{enumerate}
\qed
\end{lemma}

\begin{thm}
\label{nahmstab}
Let $\mathcal{T}=((E_1,\nabla_1),(E_2,\nabla_2),\Phi)$ be a
covariantly constant integrable triple on $C$ which satisfies the
$\tau$-coupled equations and let $(E,\nabla^{\mathcal{T}})$ be
its associated bundle with connection over $C\times \mathbb{P}^1$.
\begin{enumerate}
\item If $\tau>0$ and $\tau^\prime>0$ then the Nahm transform
$\widehat{\mathcal{T}}=((\widehat E_1, \widehat{\nabla}_1),(\widehat E_2,
\widehat{\nabla}_2),\widehat{\Phi})$ is a covariantly constant integrable
triple. Moreover, $\widehat{\mathcal{T}}$ satisfies the $\hat\tau$-coupled
equations, for some value of $\hat\tau$, if and only if $\tau=\tau^\prime$.

\item If $\tau<0$ and $\tau^\prime<0$ then the Nahm transform
$\widehat{\mathcal{T}}=((\widehat E_1, \widehat{\nabla}_1),(\widehat E_2,
\widehat{\nabla}_2),\widehat{\Phi})$ is a covariantly constant integrable
triple. Moreover, $\widehat{\mathcal{T}}$ satisfies the $\hat\tau$-coupled
equations, for some value of $\hat\tau$, if and only if $\tau=\tau^\prime$.
\end{enumerate}
\end{thm}

\begin{proof}
Proposition \ref{decomp} gives us a decomposition \begin{align*}
E_1 &\simeq\Ker \Phi^* \oplus E_1^\prime\\
E_2 &\simeq\Ker \Phi^{\ } \oplus E_2^\prime
\end{align*} Since $\mathcal{T}$ satisfies the $\tau$-couple vortex
equations, Proposition \ref{cov-const} implies that $(\Ker \Phi^*,\nabla_1)$,
$(E_1^\prime,\nabla_1)\simeq (E_2^\prime,\nabla_2)$ and $(\Ker \Phi,\nabla_2)$
are bundles with constant central curvature  with slopes $\mu(\Ker
\Phi^*)=\tau$, $\mu(E_1^\prime)=\mu(E_2^\prime)=\frac{1}{2}(\tau+\tau^\prime)$,
$\mu(\Ker \Phi) =\tau^\prime$. Now if we apply the Nahm transform and denote
$(\widehat\Phi)^*$ by $\widehat\Phi^*$, Theorem \ref{eh} implies that $({\Ker
\widehat\Phi^*}=\widehat{\Ker \Phi^*},\widehat\nabla_1)$, $(\widehat
E_1^\prime,\widehat \nabla_1)\simeq (\widehat E_2^\prime,\widehat\nabla_2)$
and $({\Ker \widehat\Phi^{}}=\widehat{\Ker \Phi^{}},\widehat\nabla_2)$ are
bundles with constant central curvature and we get a decomposition
\begin{align*}
\widehat E_1 &\simeq\Ker \widehat\Phi^* \oplus \widehat{E_1^\prime}\\
\widehat E_2 &\simeq\Ker \widehat\Phi^{\ } \oplus \widehat{E_2^\prime}
\end{align*}
The conditions $(\Phi^*\Phi)_{|E_2^\prime}=\lambda\,\mathrm{Id}_{E_2^\prime}$,
$(\Phi\Phi^*)_{|E_1^\prime}=\lambda\,\mathrm{Id}_{E_1^\prime}$  with
$\lambda\neq 0$ imply $(\widehat\Phi^*\widehat\Phi)_{|\widehat
E_2^\prime}=\lambda\,\mathrm{Id}_{\widehat E_2^\prime}$,
$(\widehat\Phi\widehat\Phi^*)_{|\widehat
E_1^\prime}=\lambda\,\mathrm{Id}_{\widehat E_1^\prime}$. Let us prove the
first equality in the $\IT_0$ case. Given $s,t\in \widehat E_{2,\xi}=\Ker
\overline\partial^{\nabla_{2,\xi}}\subset\Omega^0(E_{2,\xi})$ one has $$\langle
\widehat\Phi^*\widehat\Phi(s),t\rangle_{\widehat E_{2,\xi}}=\langle
\widehat\Phi(s),\widehat\Phi(t)\rangle_{\widehat E_{2,\xi}}$$ Taking into
account the definition of the Hermitian metric on $\widehat E_{2,\xi}$ given in
(\ref{met-herm}) of Section \ref{nahm} and the definition of $\widehat\Phi$
given in Theorem \ref{compatibility} we get $$\langle \widehat\Phi
(s),\widehat\Phi(t)\rangle_{\widehat E_{2,\xi}}=\int_{C_{\xi}}\langle
\Phi(s),\Phi(t)\rangle_{E_2}\,\omega=\int_{C_{\xi}}\langle
\Phi^*\Phi(s),t\rangle_{E_2}\,\omega$$ Therefore, if $s,t\in\widehat
E_{2,\xi}^\prime$ one has $$\langle
\widehat\Phi^*\widehat\Phi(s),t\rangle_{\widehat E_{2,\xi}}=\lambda\, \langle
s,t\rangle_{\widehat E_{2,\xi}}$$ which proves our claim. The second equality
follows in the same way. The proofs in the $\IT_1$ case are entirely similar.

This proves that $((\widehat E_1,\widehat\nabla_1), (\widehat
E_2,\widehat\nabla_2),\widehat\Phi)$ is a covariantly constant integrable
triple. Moreover, the slopes of these bundles are $\mu(\Ker
\widehat\Phi^*)=-\frac{1}{\tau}$, $\mu(\widehat E_1^\prime)=\mu(\widehat
E_2^\prime)=-\frac{2}{\tau+\tau^\prime}$, $\mu(\Ker \widehat\Phi)
=-\frac{1}{\tau^\prime}$. An easy computation shows now that
$\widehat{\mathcal{T}}$ fulfills the conditions of Proposition \ref{cov-const}
in order to have a solution of the $\hat\tau$-coupled vortex equations, for
some value of $\hat\tau$, if and only if  $\tau=\tau^\prime$.
\end{proof}

As a consequence of the preceding Theorem and the Hitchin-Kobayashi
correspondence for triples (Theorem 5.1 in \cite{BG}), which establishes an
equivalence between holomorphic triples which satisfy the $\tau$-coupled
equations and $\alpha$-polystable triples, we have.
\begin{corol}
Polystability is not preserved, in general, under the Fourier-Mukai and Nahm
transform.
\end{corol}
\begin{proof}
It is enough to take any stable bundles $F_1$, $F_2$, $F$ such that
$\mu(F)=\frac{1}{2}(\mu(F_1)+\mu(F_2))$ and $\mu(F_1)>\mu(F_2)$, which are
known to exist since the moduli spaces of stable bundles with fixed coprime
rank and degree over an elliptic curve $C$ are isomorphic to $C$ and thus they
are not empty (see \cite{Tu}). Now define the triple $T=(F_1,0,0)\oplus
(F,F,\mathrm{Id}_F)\oplus (0,F_2,0)$ and endow $F_1$, $F_2$, $F$ with
connections of constant curvature compatible with their holomorphic structures
according to Donaldson Theorem \cite{Do}. Now, Proposition \ref{cov-const}
implies that $T$ is $(\tau-\tau^\prime)$-polystable since, by construction, it
satisfies the $\tau$-coupled equations, with $\tau=\mu(F_1)$ and
$\tau^\prime=\mu(F_2)$.

If we take $\mu(F_1)\neq \mu(F_2)$,  Theorem \ref{nahmstab} implies that the
transformed triple $\widehat T$ does not satisfy the $\hat\tau$-coupled
equations for any value of $\hat\tau$. By the Hitchin-Kobayashi
correspondence for triples \cite[Theorem 5.1]{BG}, this implies that
$\widehat T$ is not polystable.
\end{proof}

The preservation of stability remains as an open question. Notice that in the
case of stable triples $(E_1,E_2,\Phi)$ with $E_1\neq 0$ and $E_2\neq 0$, the
condition of being covariantly constant implies that $\Phi$ is an isomorphism
(Remark  \ref{isocov}). Now stability is preserved in the conditions of Theorem
\ref{suinvone}.

{\bf Acknowledgement.} The authors are deeply grateful to the anonymous referee
for his/her remarks which have helped to significantly improve the paper.

\end{document}